%% file: main.tex
\documentclass[authoryear,preprint,review,11pt]{elsarticle}

\journal{European Journal of Operational Research}
\biboptions{authoryear}
\usepackage{hyperref}
\hypersetup{
  colorlinks   = true, 
  urlcolor     = blue, 
  linkcolor    = blue, 
  citecolor   = blue!70 
}

\usepackage{setspace}

\usepackage{graphicx}

\usepackage{float}
\usepackage{multirow}
\usepackage{subfigure}
\usepackage{verbatim}
\usepackage[dvipsnames]{xcolor}
\usepackage[figuresright]{rotating}
\usepackage[normalem]{ulem}
\usepackage{appendix}

\usepackage{url}
\usepackage{enumitem}

\graphicspath{ {figures/} }

\usepackage{relsize}

\usepackage[utf8]{inputenc}
\usepackage{amsmath,amssymb,amsthm}
\usepackage{bm}

\usepackage{pdflscape}
\usepackage{afterpage}
\usepackage{capt-of}
\usepackage{fullpage}
\usepackage{tikz}

\usepackage[ruled,commentsnumbered,linesnumbered,vlined]{algorithm2e}
\SetKw{KwAnd}{and}
\SetKw{KwOr}{or}
\SetKw{KwNot}{not}
\SetKw{KwContinue}{continue}
\SetKw{KwTrue}{true}
\SetKw{KwFalse}{false}
\SetKw{KwNull}{null}
\SetKw{KwMod}{mod}
\SetKwRepeat{Do}{do}{while}

\usepackage{array}     

\usepackage{placeins}

\usepackage{indentfirst}
\usepackage[normalem]{ulem}
\usepackage{textcomp}

\usepackage{wasysym} 
\usepackage{geometry}
\usepackage{soul}

\usepackage{booktabs, multirow} 
\usepackage{xcolor,colortbl} 
\usepackage{changepage,threeparttable} 

\usepackage{interval}
\usepackage{media9}

\usepackage{mathtools}
\DeclarePairedDelimiter\ceil{\lceil}{\rceil}
\DeclarePairedDelimiter\floor{\lfloor}{\rfloor}

\definecolor{napiergreenMod}{rgb}{0.16, 0.6, 0.0}
\definecolor{myOrange}{cmyk}{0, 0.533, 1, 0.15}

\newcommand{\vitorC}[1]{\textcolor{black}{#1}}

\newcommand{\rafaelF}[1]{\footnote{\textcolor{red}{{\sc rafael}: #1}}}

\newcommand{\okC}[1]{\textcolor{black}{#1}}

\let\citeA\cite        
\let\cite\citep        

\usepackage{hyphenat} 
       \makeatletter
\def\@fnsymbol#1{\ensuremath{\ifcase#1\or *\or \mathsection\or \mathparagraph\or **\or \dagger\or \ddagger\or
    \|\or \dagger\dagger
   \or \ddagger\ddagger \else\@ctrerr\fi}}
    \makeatother

\date{\today}

\begin{document}

\begin{frontmatter}

\title{Constraint programming model and biased random-key genetic algorithm for the single-machine coupled task scheduling problem with exact delays to minimize the makespan}

\author[1st_address]{Vítor A. Barbosa}
\ead{vitor.barbosa@ufba.br}

\author[1st_address]{Rafael A. Melo\corref{cor1}}\cortext[cor1]{Corresponding author}
\ead{rafael.melo@ufba.br}

\address[1st_address]{Institute of Computing, Universidade Federal da Bahia, Salvador, BA 40170-115, Brazil}

\begin{abstract}
We consider the \okC{strongly NP-hard} single-machine coupled task scheduling problem with exact delays to minimize the makespan. 
In this problem, a set of jobs has to be scheduled, each composed of two tasks interspersed by an exact delay. 
\okC{Given that no preemption is allowed, the goal consists of minimizing the completion time of the last scheduled task.} 
We model the problem using constraint programming (CP) and propose a biased random-key genetic algorithm (BRKGA).
\okC{Our CP model applies well-established global constraints.}
Our BRKGA combines some successful components in the literature: 
\okC{an initial solution generator}, periodical restarts and shakes, and a local search algorithm.
\okC{Furthermore, the BRKGA's decoder is focused on efficiency rather than optimality, which accelerates the solution space exploration.}
Computational experiments \okC{on a benchmark set containing instances with up to 100 jobs (200 tasks)} indicate that the proposed BRKGA can \okC{efficiently} explore the problem solution space\okC{, providing high-quality approximate solutions within low computational times}. 
It can also provide better solutions than the CP model under the same computational settings, i.e., three minutes of time limit and a single thread.
The CP model, when offered a longer running time of 3600 seconds and multiple threads, significantly improved the results, \okC{reaching the current best-known solution for 90.56\% of these instances.}
Finally, our experiments highlight the importance of the shake and local search components in the BRKGA, whose combination significantly improves the results of a standard BRKGA. 
\end{abstract}

\end{frontmatter}

\setlength{\unitlength}{4144sp}%
\begingroup\makeatletter\ifx\SetFigFont\undefined%
\gdef\SetFigFont#1#2#3#4#5{%
  \reset@font\fontsize{#1}{#2pt}%
  \fontfamily{#3}\fontseries{#4}\fontshape{#5}%
  \selectfont}%
\fi\endgroup%

\input{050_introduction}
\input{100_constraintprogramming}

\input{150_brkga}

\input{160_brkga_single}
\input{165_brkga_notations_data_structures_decoder}
\input{170_brkga_decoder_single}
\input{173_brkga_single_decoder_first_fit_analysis}

\input{175_brkga_single_warm_start}
\input{180_brkga_single_local_search}

\input{182_brkga_ls_neighborhood}

\input{183_brkga_ls_algorithm}

\input{184_brkga_ls_integration}
\input{185_brkga_single_perturbation}
\input{200_brkga_single_overall_framework}

\input{250_computationalresults}

\input{255_benchmark_instances}

\input{260_parametrization}

\input{263_comparison_cp_vs_mips_local}
\input{272_comparison_our_approaches_vs_bks}

\input{265_comparison_tested_approaches}

\input{283_impact_of_the_brkga_components}
\input{300_finalcomments}

\section*{Acknowledgments}
This study was financed in part by the Coordenação de Aperfeiçoamento de Pessoal de Nível Superior - Brasil (CAPES) - Finance Code 001.
This work was partially supported by the Brazilian National Council for Scientific and Technological Development (CNPq) research \okC{grants} 314718/2023-0 \okC{and 445324/2024-4}, and the FAPESB INCITE PIE0002/2022 grant.

\bibliography{main}

 \newpage

 \begin{appendices}
 \input{800_appendix_irace}

 \input{810_population_homogeneity}

 \end{appendices}

\end{document}

%% file: 050_introduction.tex
\section{Introduction}
\label{sec:intro}

The coupled-task scheduling problem (CTSP) consists of scheduling jobs, each with at least two coupled tasks.
In each job, the succeeding task must be started after the preceding task is finished, with an exact delay between them.
According to \citeA{KhaSalChe20}, there are several possible performance criteria for the CTSP.
However, the major performance criteria adopted in the literature include the makespan, total completion time, maximum lateness, and earliness and tardiness.
For this work, we tackle the single-machine version of the problem, where each job has two tasks that must be separated by an exact delay.
The goal consists of minimizing the makespan.
\okC{Such minimization is challenging because tasks can be scheduled within the exact delay of a job, whether by nesting jobs between other job tasks or by interleaving tasks of different jobs.
Henceforth, there must be several non-singleton jobs, i.e., jobs with exact delays larger than their processing times.
Otherwise, these jobs can be randomly appended at the end of a schedule without loss of generality \cite{LiZha07}.}

Applications of the CTSP include pulsed radar systems in which a radar transmits a pulse and receives its reflection to keep track of targets \cite{Sha80}, chemistry manufacturing processes~\cite{AgeBab07}, patient appointments in medical services~\cite{Peretal11, ConSha14, Liuetal19}, among others.
Another potential application involves the scheduling of tasks operated by robots in smart homes, especially in smart kitchens \cite{Sharma2024}. 
In this scenario, kitchen robots~\cite{MEPANI20221930, Jiang2022} can take advantage of idle moments between sub-tasks to perform other sub-tasks, thus reducing the total preparation time \cite{Bautista2023, Yi2022}. 
\citeA{Yi2022} studied a similar problem involving multiple tasks operated by dual-arm kitchen robots in a controlled environment.

The single machine coupled-task scheduling problem (SMCTSP) to minimize makespan is classified as a \okC{strongly} NP-hard problem \cite{OrmPot97, SheSmi05}. 
\citeA{ConSha12} showed that obtaining an optimal solution is NP-hard even if there is a predefined sequence for the initial (or final) tasks and all the jobs have unit processing times.
\okC{\citeA{LiZha07} proposed a tabu-search heuristic for the SMCTSP.
\citeA{HWANG2011690} studied the SMCTSP with fixed-job-sequence, i.e., if job $A$ precedes job $B$, then the initial and final tasks of job $A$ must precede the initial and final tasks of job $B$, respectively.
The authors designed an $O(n^2)$ algorithm to schedule with minimum makespan a given task sequence abiding by the fixed-job-sequence constraint, in addition to determining its feasibility.}
\citeA{BesGir19} recently investigated the parameterized complexity of the CTSP.
\citeA{KhaSalChe20} presented a comprehensive review of the different types of CTSP.
They also presented mixed integer programming (MIP) formulations for several coupled scheduling problems.
\citeA{KhaSal21} presented a binary search-based MIP heuristic\okC{, in addition to some lower and upper bounds}.
\citeA{BekDosGal22} provided an approximation algorithm for a special case in which jobs only have two different delay times.
Recently, \citeA{KhaSal24} proposed a new IP formulation and a relax-and-solve heuristic for the SMCTSP, which, \okC{to the best of our knowledge}, is the current state-of-the-art method.
\okC{We refer to \citeA{KhaSalChe20} for a literature survey on coupled task scheduling problems.}

The main contributions of our work are highlighted in the following.
First, we propose a constraint programming (CP) model that outperforms all the state-of-the-art approaches, both exact and heuristic methods.
Second, we devise a BRKGA metaheuristic, denoted as BRKGA-R-S-LS, that can efficiently explore the solution space and find high-quality solutions in low computational times. 
Besides, it outperforms the CP model when both are run under the same computational settings, i.e., a low time limit and a single thread.
\okC{It is worth highlighting that we considered the 240  instances available from the literature, which have instances with up to 100 jobs (200 coupled tasks) \cite{csbenchmark, KhaSalChe20}.
Over them, the CP model obtained 71 optimal solutions.
Considering the solutions obtained in this paper, there are 73 optimal solutions in total.}

Our BRKGA-R-S-LS stands for \textit{Biased Random-Key Genetic Algorithm with Restarts, Shakes and Local Search}, and combines some successful components in the literature: 
\okC{an initial solution generator}, periodical restarts and shakes, and a local search procedure.
The \okC{initial solution generator} aims to provide the competitive advantage of starting the search guided by good solutions.
The restarts and shakes intend to direct the search to new paths, whether completely different (restart) or similar (shake).
Finally, the objective of the local search is to assist the intensification mechanism of the search.
The BRKGA-R-S-LS is also compared to a standard BRKGA with \okC{an initial solution generator} restarts, denoted as BRKGA-R, to analyze the impact of the shake and local search on solution quality.

\subsection{\okC{Constraint programming motivation}}

\okC{CP is a computational paradigm for solving a constraint satisfaction problem (CSP)~\cite{rossi2006}. A CSP is represented in terms of decision variables and the relations that must hold among them, i.e., the constraints.
A constraint solver takes one CSP instance and assigns a value for each variable while satisfying the problem constraints.
When a feasible solution is achieved, the solver propagates the constraints to prune the search space, reducing the domain of the variables.
Besides, each constraint solver applies dedicated algorithms for exploiting the structure and properties of the constraint.}

\okC{Note that CP is focused on feasibility rather than optimality.
Nevertheless, CP can also be used to find an optimal solution for optimization CSPs by solving a sequence of CSPs.
Therefore, the constraint solver will initially find a feasible solution for the problem with a value.
Then, it imposes a constraint that the next solution must be better.
Whenever the solver fails to find a feasible solution, the previous solution is optimal.
Although CP can take too long in practice, it is guaranteed to find an optimal solution.
Additionally, it has been demonstrated to be a strong approach for scheduling problems, because the new constraint reduces the search space.}

\okC{The reason to apply CP to the SMCTSP is its potential success over scheduling problems in the past years. 
For instance, \citeA{LUNARDI2020105020} devised CP models for the online printing shop scheduling problem, solving to optimality all small-sized and a fraction of the medium-sized instances, in addition to providing feasible solutions for large-sized instances that appear in practice.
\citeA{DACOL2022100249} showed that CP can also be efficiently applied to industrial-size job scheduling problems, dealing with instances with up to one million operations.
\citeA{AWAD2022107565} demonstrated the success of their CP model to minimize the makespan of a large scheduling problem in batch pharmaceutical manufacturing facilities.}

\subsection{\okC{Biased random-key genetic algorithm motivation}}

\okC{The BRKGA is a population-based search metaheuristic proposed by \citeA{gonccalves2011biased} that simplifies genetic algorithms by making the solution representation aspects and the intensification-diversification mechanism problem independent.
Such simplification led the BRKGA to be successfully applied to several combinatorial optimization problems.
Besides, the structure supports hybridization with other optimization methods in different scenarios.}

\okC{Recently, \citeA{melo2023biased} and \citeA{silva2024obtaining} devised BRKGAs to obtain high-quality solutions in low computational times to the minimum quasi-clique partitioning problem and the Grundy coloring problem, respectively.
\citeA{Reixach2025} proposed a BRKGA hybridized with a beam search from the literature to the longest common square subsequence problem, outperforming former state-of-the-art methods with statistical significance.
\citeA{Chaves2024} developed an adaptive BRKGA for the tactical berth allocation problem, which was competitive with the current state-of-the-art methods and could find the best-known solutions for most tested instances.
Concerning scheduling problems, nearly one-third of all BRKGA applications are dedicated to this category of problems, making it the most extensively studied category using this metaheuristic \citep{LONDE20251}.}

\okC{In the context of hybridization, one or multiple heuristics can be applied to BRKGA to provide one or multiple \okC{initial} solutions.
This type of solution is included in the initial population and must meet some quality criterion.
\citeA{MARTARELLI2020100618} and \citeA{samuelCelsoUeverton2023} showed that this strategy indicated a positive influence on the performance of the BRKGA, given the competitive advantage of starting the search guided by good solutions. 
Another hybridization involves employing a local search heuristic, which performs consecutive minor changes on a solution.
These minor changes, called \textit{moves}, are limited by the type of operation defined a priori, i.e., there is a constrained set of solutions to be reached, which we call \textit{neighborhood}.
Given that each move must improve the current solution, the search finishes when no possible move results in an improvement.
At this point, the solution is called a local optimum.
The local search heuristic can be applied in several scenarios.
For instance, \citeA{silva2024} executed the local search after an improvement in some individuals of the elite set, while \citeA{DEABREU2022101149} applied the local search only on selected iterations.}

\okC{Thereupon, the motivation behind the application of the BRKGA to the SMCTSP relies on its success over scheduling problems.
More specifically, the single machine coupled-task scheduling problem is constrained, considering that it imposes an exact delay between the tasks of the jobs, and the BRKGA can always decode feasible solutions with low computational effort.
Furthermore, once a solution is decoded, no manipulation is necessary since the BRKGA will independently deal with the intensification and diversification mechanisms.
Finally, its hybridization with the local search generally yields a positive impact on scheduling problems. 
We refer the reader to \citeA{Londe2024} and \citeA{LONDE20251} for recent literature surveys about the BRKGA.}

\subsection{Organization}

The remainder of this paper is organized as follows.
Section~\ref{sec:notationscp} presents the notations used and the CP model.
Section~\ref{sec:brkga} describes the components of the proposed BRKGA-R-S-LS metaheuristic.
Section~\ref{sec:experiments} reports the results obtained from a series of computational experiments on the proposed methods.
Section~\ref{sec:concluding} concludes the paper with some final remarks.

%% file: 100_constraintprogramming.tex
\section{Notations and constraint programming model}
\label{sec:notationscp}

\begin{figure}[h]
    \centering
    \includegraphics[width=.9\linewidth, trim={0.8cm 0 1.28cm 0}, clip]{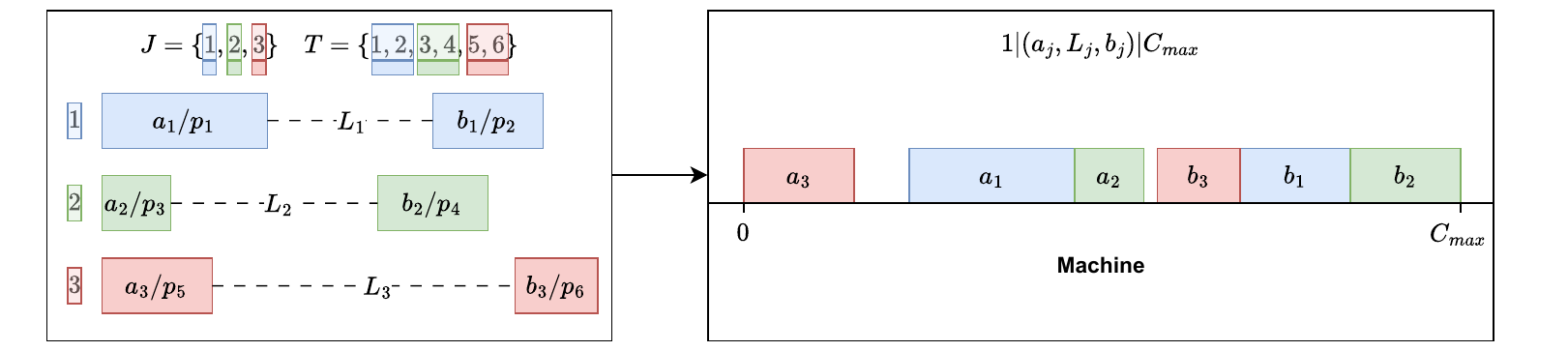}
    \caption{Illustration of the single machine coupled task scheduling problem}
    \label{fig:smctsp_illustration}
\end{figure}

Consider a set of jobs $J=\{1,\ldots,n\}$. 
Each job $j \in J$ is represented by $(a_j, L_j,b_j)$, where $a_j \in \mathbb{Z}_{\geq 0}$ and $b_j \in \mathbb{Z}_{\geq 0}$ denote the processing times of their initial and final tasks, respectively. 
$L_j \in \mathbb{Z}_{\geq 0}$ provides the exact interval between the two tasks of job $j$.
Preemption is not allowed, and the objective is to minimize the makespan $C_{max}$, i.e., the completion time of the last scheduled job.
\okC{Considering standard scheduling classifications, the problem can be denoted as $1|(a_j, L_j,b_j)|C_{max}$, where $1$ indicates a single machine, $(a_j, L_j,b_j)$ provides the jobs' characteristics, and $C_{max}$ defines the objective function}.
It is assumed that processing times and intervals between tasks are integer values.
Additionally, denote the set of tasks by $T = \{1,\ldots,2n\}$, where tasks $2j-1$ and $2j$ identify the initial and final tasks of job $j \in J$, respectively. 
Let $p_h$ denote the processing time of task $h \in T$.
Figure \ref{fig:smctsp_illustration} illustrates the single machine coupled task scheduling problem along with its notations.

Define the variable $s_h$ to represent the start time of task $h \in T$.
Let the variable $C_{max}$ indicate the completion time of the last scheduled task.
Denote by $UB$ an upper bound for $C_{max}$, which a possibility is to define as $UB = \sum_{j \in J}(a_j + L_j + b_j)$. 
The problem can be formulated as the following \okC{CP} model:
\begingroup
\allowdisplaybreaks
\begin{align}
 & \min \  C_{max}  \label{cp:obj}\\
 & s_{2j} = s_{2j-1} + p_{2j-1} + L_j, \qquad \forall \ j \in J,       \label{cp:01} \\
& \textrm{noOverlap}(\{s_h | h \in T\},\ \{p_h | h \in T\}),    \label{cp:03} \\
& C_{max} = \max_{h \in T} (s_{h} + p_{h}),   \label{cp:04} \\
&  s_h \in [0, UB], \qquad \forall \ h \in T, \label{cp:05}\\
&  C_{max} \in [0, UB]. \label{cp:07}
\end{align}
\endgroup


The objective function~\eqref{cp:obj} minimizes the completion time of the last scheduled task, i.e., the makespan.
Constraints~\eqref{cp:01} ensure the exact interval between the two tasks of each job. 
\okC{The global constraint \mbox{noOverlap(A, B)} in Constraint~\eqref{cp:03} ensures that the tasks in $T$ do not overlap, in which the starting times $A=\{s_h | h \in T\}$ and their respective processing times $B=\{p_h | h \in T\}$ represent the intervals $[s_h, s_h+p_h], \ h \in T,$ that must be disjoint.
It is noteworthy that the method  \mbox{noOverlap(A, B)} can be initially modeled as a conjunction of non-overlap constraints, but can also be modeled differently depending on the underlying solver.} Constraint~\eqref{cp:04} defines the makespan as the completion time of the last task. Constraints~\eqref{cp:05}-\eqref{cp:07} define the domains of the decision variables as integer values in the specified intervals.

%% file: 150_brkga.tex
\section{Biased random-key genetic algorithm}
\label{sec:brkga}

\okC{The BRKGA is based on the concept of Genetic Algorithms (GAs), where each problem's solution is an individual within a population, and its solution value is denoted as \textit{fitness}.
Based on the random-key genetic algorithm (RKGA) of \citeA{bean1994genetic}, the solutions are encoded into vectors of real numbers generated randomly in the interval $[0,1)$, these numbers being denoted keys. 
In each iteration of the BRKGA, the current population generation is decoded using a deterministic algorithm, referred to as decoder, sorted by fitness, and partitioned into two subsets: elite and non-elite.
The elite set contains the individuals with the best fitness, while the non-elite set comprises the remaining solutions.}

The intensification mechanism in a BRKGA comprises the crossover process within the creation of a new generation of solutions.
In this process, one of the parents is always an elite solution, defined as one with a high fitness value. 
In our case, the fitness value can be defined as a negative makespan.
Moreover, the elite parent has a greater probability of passing on its characteristics, i.e., keys, to the new population.
The BRKGA uses the parametric uniform crossover strategy~\cite{Spe91} to combine parent solutions and generate a new one.

For diversification, with each new iteration, the BRKGA adds a set of randomly generated solutions, called mutants, to the population.
Furthermore, a BRKGA approach can allow restarts, i.e., a new initial population is generated from which the algorithm must continue generating solutions.
The BRKGA can also incorporate a shake operation \cite{ANDRADE201967}, which perturbs all individuals from the elite set and resets the remaining population.

In this work, we propose a BRKGA that combines some successful methods in the literature:
\okC{an initial solution generator}, periodical restarts or shakes, and a local search. It is denoted henceforth as BRKGA-R-S-LS.
To highlight the importance of a shaking operation combined with a local search, we compare it with a standard BRKGA with \okC{an initial solution generator} and restarts, which we denote as BRKGA-R.

Throughout this section, we describe each component used in BRKGA-R-S-LS, namely: 
the encoding and the decoder (Section \ref{sec:encoding_decoder_brkga}); 
the \okC{initial solution generator} (Section \ref{sec:warm-start}); 
the local search (Section \ref{sec:local_search}); 
and the perturbation (Section \ref{sec:perturbation}).
\okC{At the end of the section, we detail the overall framework (Section~\ref{sec:overall_framework}).}

%% file: 160_brkga_single.tex
\vspace{-0.3cm}
\subsection{Encoding and decoder}
\label{sec:encoding_decoder_brkga}

In our approach, each solution is encoded as a vector $X$ composed of $|J|$ random keys, in which the \okC{$j$-th key corresponds to the job $j \in J$}.

\begin{figure}[!ht]
    \centering
    \includegraphics[width=.85\linewidth, trim={0.8cm 0.1cm 1.28cm 0.1cm}, clip]{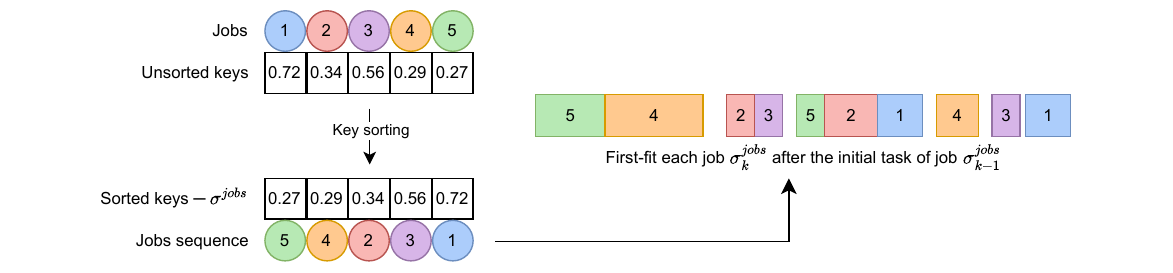}
    \caption{Decoder example.}
    \label{fig:sorted_sequence}
\end{figure}

Our decoder is a polynomial-time heuristic in which the keys are used to define a sorted sequence of the jobs $\sigma^{jobs} = (\sigma^{jobs}_1, \sigma^{jobs}_2, \ldots, \sigma^{jobs}_n)$, where $\sigma^{jobs}_k, \ k \in \{1,\ldots,n\}$, is the $k$-th job in the sequence. 
\okC{The job $\sigma^{jobs}_k$ is thus inserted into the solution after the finish of the initial task of $\sigma^{jobs}_{k-1}$ in the earliest possible time without overlapping with those jobs already scheduled (first-fit position).
By construction, our first-fit algorithm is injective, i.e., given two different sorted sequences of the jobs, the resulting solutions are also different.
On the other hand, our decoder is not injective because two random-key vectors can produce the same sorted sequence of the jobs.
Figure \ref{fig:sorted_sequence} illustrates a decoding.}

\FloatBarrier
\vspace{-0.5cm}

%% file: 165_brkga_notations_data_structures_decoder.tex
\subsubsection{Notations and main data structures for the decoder}
\label{sec:notations_and_data_structures}

The decoder algorithm uses the following notations and data structures, which are summarized in Table \ref{tab:summary_notation_decoder}.
Let the tuple $S: \langle \tau, \sigma^{tasks}, C_{max} \rangle$ represent a solution for the SMCTSP.
Each solution $S$ contains the following structure: 
$S.\tau$ is an array of integers of length $|T|$ in which $S.\tau_h$ defines the start time of the task $h \in T$ in $S$;
$S.\sigma^{tasks}$ is an array of integers of length $|T|$ in which $S.\sigma^{tasks}_i$, $i \in \{1,\ldots,|S.\sigma^{tasks}|\}$ is the $i$-th task of the schedule in $S$;
and $S.C_{max}$ is the makespan of $S$.

Our decoder iteratively inserts pairs of tasks (jobs) in the schedule $S.\sigma^{tasks}$.
Therefore, each pair of positions to insert a job is called an insertion candidate.
This insertion candidate has a specific structure, further used to update $S$.
Let the tuple $\mathcal{C}: \langle j,\, F,\, pos_{t_1}, \, pos_{t_2}, \, st, \, cost \rangle$ be an insertion candidate for a specific job.
Each insertion candidate $\mathcal{C}$ for a job $j \in J$ contains the following structure:
\okC{$\mathcal{C}.j$} is the id of the job to be inserted;
$\mathcal{C}.F$ is a boolean variable that indicates its feasibility;
$\mathcal{C}.pos_{t_1}$ and $\mathcal{C}.pos_{t_2}$ define the insertion position of the initial and final tasks of \okC{$\mathcal{C}.j$}, respectively;
$\mathcal{C}.st$ defines $j$'s start time;
$\mathcal{C}.cost$ defines its cost, which is given by the difference between the new and current values of $S.C_{max}$.
\okC{It is worth mentioning that when $\mathcal{C}.pos_{t_1} = \mathcal{C}.pos_{t_2}$, the whole job $\mathcal{C}.j$ is inserted between the tasks $\mathcal{C}.\sigma^{tasks}_{pos_{t_1}-1}$ and $\mathcal{C}.\sigma^{tasks}_{pos_{t_1}}$.
Additionally, if $\mathcal{C}.pos_{t_h} = |S.\sigma^{tasks}|+1$, then the respective task $h \in \{1,2\}$ (initial or final task) is inserted at the end of the schedule $S.\sigma^{tasks}$.}

\vspace{-0.5cm}

\begin{table}[!htp]\centering
\caption{Summary of the notation used in the decoder} \label{tab:summary_notation_decoder}
\setstretch{1.0}
\scriptsize
\begin{tabular}{c|l|l|l}\hline
\multicolumn{1}{l|}{Structure}& \multicolumn{1}{l|}{Definition}&\multicolumn{1}{c|}{Notation}  & \multicolumn{1}{l}{Definition} \\\hline
\multirow{3}{*}{$S$} &\multirow{3}{*}{solution} &$S.\tau$ & array of the task start times of $S$\\
&&$S.\sigma^{tasks}$ & array of the sequence of tasks in $S$\\
&&$S.C_{max}$ &makespan of $S$ \\
\hline
\multirow{6}{*}{$\mathcal{C}$} & \multirow{6}{*}{insertion candidate} 
&$\mathcal{C}.j$ & job id\\
&&$\mathcal{C}.F$ & true if $\mathcal{C}$ is feasible, false otherwise\\
&&$\mathcal{C}.pos_{t_1}$ &insert position for the initial task \\
&&$\mathcal{C}.pos_{t_2}$ &insert position for the final task \\
&&$\mathcal{C}.st$ &start time of the job $\mathcal{C}.j$\\
&&$\mathcal{C}.cost$ &cost of $\mathcal{C}$ \\
\hline
\end{tabular}
\end{table}

%% file: 170_brkga_decoder_single.tex
\subsubsection{First-fit algorithm}

Algorithm~\ref{alg:single_dec} details the first-fit algorithm.
It receives as input the insertion sequence of the jobs $\sigma^{jobs}$ and the instance data $I$. 
The instance data $I$ contains the parameters presented in Section~\ref{sec:notationscp}.
In line~\ref{alg:single_dec:010}, a partial solution~$S$ is initiated with the job $\sigma^{jobs}_1$ scheduled at time 0, i.e., as early as possible.
\okC{In line~\ref{alg:single_dec:020}, we define a variable $\overleftarrow{pos}_{t_1}$ that stores the position of the last initial task inserted in $S.\sigma^{tasks}$, i.e., $\overleftarrow{pos}_{t_1} = 1$ because the initial task of job $\sigma^{jobs}_1$ was inserted in the first position of $S.\sigma^{tasks}$ (line~\ref{alg:single_dec:010}).
In the for loop of lines~\ref{alg:single_dec:030}-\ref{alg:single_dec:120}, we iteratively insert the remaining jobs of the sequence $\sigma^{jobs}$ (line~\ref{alg:single_dec:030}) into the schedule of the solution~$S$.
In line~\ref{alg:single_dec:035}, we define the job to be inserted as the job $\sigma^{jobs}_i$.
In line~\ref{alg:single_dec:040}, we define a dummy infeasible insertion candidate $\mathcal{C}^*$ to serve as a sentinel.
In line~\ref{alg:single_dec:050}, we define another variable $pos_{t_1}$ that stores the expected position of the initial task of job $\sigma^{jobs}_i$.
By definition, $pos_{t_1} \geq \overleftarrow{pos}_{t_1} + 1$ because the job $\sigma^{jobs}_i$ must be inserted after the finish of the initial task of job $\sigma^{jobs}_{i-1}$.
In the loop of lines~\ref{alg:single_dec:060}-\ref{alg:single_dec:080}, we iteratively search for the first-fit insertion candidate $\mathcal{C}^*$ when inserting the initial task of job $\sigma^{jobs}_i$ at position $pos_{t_1} \leq |S.\sigma^{tasks}|$ (line~\ref{alg:single_dec:060}).
In line~\ref{alg:single_dec:070}, the function \mbox{ANALYZE\_INSERTION\_CANDIDATE}, described in Algorithm~\ref{alg:analyze_candidate}, analyzes the insertion candidate's feasibility, cost, and final task's position, assuming that the initial task is placed in position $pos_{t_1}$.
In line~\ref{alg:single_dec:080}, we increment the variable $pos_{t_1}$.
At the end of the loop of lines~\ref{alg:single_dec:060}-\ref{alg:single_dec:080}, the variable $\mathcal{C}^*$ can be either an infeasible or a feasible insertion candidate.
If it is infeasible (line~\ref{alg:single_dec:090}), the job $\sigma^{jobs}_i$ can only be inserted in the position $pos_{t_1} = |S.\sigma^{tasks}|+1$ with cost $a_j + L_j + b_j$ (line~\ref{alg:single_dec:100}).
After that, the variable $\mathcal{C}^*$ is a feasible insertion candidate and we can update the solution~$S$.
In lines~\ref{alg:single_dec:103} and~\ref{alg:single_dec:106}, we store the start times of the initial and final tasks of the job $j$ in the array $S.\tau$, respectively.
In lines~\ref{alg:single_dec:109} and~\ref{alg:single_dec:112}, we execute a function $\texttt{insert(A,P,E)}$, which inserts the integer \texttt{E} in the position \texttt{P} of the array \texttt{A}.
Given that both $\mathcal{C}^*.pos_{t_1}$ and $\mathcal{C}^*.pos_{t_2}$ are positions in the currently unmodified schedule $S.\sigma^{tasks}$, the order of insertion matters.
Consequently, we first insert the final task (line~\ref{alg:single_dec:109}), then the initial task (line~\ref{alg:single_dec:112}).
In line~\ref{alg:single_dec:115}, we update the solution makespan $S.C_{max}$ by adding the cost $\mathcal{C}^*.cost$.
In line~\ref{alg:single_dec:120}, the variable $\overleftarrow{pos}_{t_1}$ is updated to $\mathcal{C}^*.pos_{t_1}$.}
Finally, the algorithm returns the solution~$S$ in line~\ref{alg:single_dec:130}.

\begin{algorithm}[H]
    \scriptsize
    \setstretch{0.75}
    \caption{\scriptsize FIRST\_FIT\_ALGORITHM($\sigma^{jobs}$, $I$)}
    \label{alg:single_dec}
    \okC{
    $S \gets$ Partial solution with job $\sigma^{jobs}_1$ scheduled at time 0\; \label{alg:single_dec:010}
    $\overleftarrow{pos}_{t_1} \gets$ 1\; \label{alg:single_dec:020}
    \For{$i \gets 2$ \KwTo $n$}{\label{alg:single_dec:030}
        $j \gets \sigma^{jobs}_{i}$\;\label{alg:single_dec:035}
        $\mathcal{C}^* \gets \mathcal{C}: \langle j, \, \KwFalse, \, -, \, -, \, -, \, - \rangle$\; \label{alg:single_dec:040}
        $pos_{t_1} \gets \overleftarrow{pos}_{t_1} + 1$\; \label{alg:single_dec:050}
        \While{
            $\mathcal{C}^*.F =$ \KwFalse \KwAnd $pos_{t_1} \leq |S.\sigma^{tasks}|$
        }{ \label{alg:single_dec:060}
            $\mathcal{C}^* \gets$ ANALYZE\_INSERTION\_CANDIDATE($j$, $pos_{t_1}$, $S$, $I$)\; \label{alg:single_dec:070}
            $pos_{t_1} \gets pos_{t_1} + 1$\; \label{alg:single_dec:080}
        }
        \If{$\mathcal{C}^*.F =$ \KwFalse}{ \label{alg:single_dec:090}
            $\mathcal{C}^* \gets \mathcal{C}: \langle j, \, \KwTrue, \, pos_{t_1}, \, pos_{t_1}, \, S.C_{max}, \, a_j + L_j + b_j \rangle$\; \label{alg:single_dec:100}
        }
        $S.\tau_{2j-1} \gets \mathcal{C}^*.st$\; \label{alg:single_dec:103}
        $S.\tau_{2j} \gets \mathcal{C}^*.st + a_j + L_j$\; \label{alg:single_dec:106}
        \texttt{insert(}$S.\sigma^{tasks}$, $\mathcal{C}^*.pos_{t_2}$, $2j$\texttt{)}\; \label{alg:single_dec:109}
        \texttt{insert(}$S.\sigma^{tasks}$, $\mathcal{C}^*.pos_{t_1}$, $2j-1$\texttt{)}\; \label{alg:single_dec:112}
        $S.C_{max} \gets S.C_{max} + \mathcal{C}^*.cost$\; \label{alg:single_dec:115}
        $\overleftarrow{pos}_{t_1} \gets \mathcal{C}^*.pos_{t_1}$\; \label{alg:single_dec:120}
    }
    \Return $S$\; \label{alg:single_dec:130}
    }
\end{algorithm}

Algorithm~\ref{alg:analyze_candidate} describes the pseudocode for the ANALYZE\_INSERTION\_CANDIDATE invoked in line~\ref{alg:single_dec:070} of Algorithm~\ref{alg:single_dec}. 
It analyzes the cost of inserting a job $j$, whose initial task is in position $pos_{t_1}$ of  $S.\sigma^{tasks}$.
In line~\ref{alg:analyze_candidate:010}, we define $prev$ as the id of the task preceding the task at position $pos_{t_1}$ in the schedule $S.\sigma^{tasks}$.
Next, in line~\ref{alg:analyze_candidate:020} we define $st_{t_1}$  as the start time of the initial task of $j$, which is right after the conclusion of task $prev$.
If the initial task of $j$ overlaps $S.\sigma^{tasks}_{pos_{t_1}}$ when starting at $st_{t_1}$ (line~\ref{alg:analyze_candidate:030}), then we return an infeasible insertion candidate with infinity cost (line~\ref{alg:analyze_candidate:040}).
In line~\ref{alg:analyze_candidate:050}, we calculate $st_{t_2}$, which is the start time of $j$'s final task.

The following lines (\ref{alg:analyze_candidate:050}-\ref{alg:analyze_candidate:160}) of Algorithm~\ref{alg:analyze_candidate} will search for the position of $j$'s final task in $S.\sigma^{tasks}$.
Therefore, if $j$'s final task fits before $S.\sigma^{tasks}_{pos_{t_2}}, \ pos_{t_2} = pos_{t_1}$ (line~\ref{alg:analyze_candidate:060}), when starting at $st_{t_2}$ (line~\ref{alg:analyze_candidate:070}), \vitorC{then} this insertion candidate is feasible with zero cost.
Otherwise, if $st_{t_2}$ is greater than the makespan of $S$ (line~\ref{alg:analyze_candidate:090}), we return a feasible insertion candidate (line~\ref{alg:analyze_candidate:120}), whose final task is scheduled at the end of $S.\sigma^{tasks}$ (line~\ref{alg:analyze_candidate:100}) with cost $st_{t_2} + b_j - S.C_{max}$ (line~\ref{alg:analyze_candidate:110}).
If none of the two possibilities holds, it means that $j$'s final task could only fit in the interval [$pos_{t_1}+1$, $|S.\sigma^{tasks}|$].
Thus, in line~\ref{alg:analyze_candidate:130}, we iteratively search for the expected position $h$ of $j$'s final task in $S.\sigma^{tasks}$, which is in the interval $[pos_{t_1}+1, |S.\sigma^{tasks}|]$.
If $j$'s final task fits before $S.\sigma^{tasks}_h$ when starting at $st_{t_2}$ (line~\ref{alg:analyze_candidate:140}), we return a feasible insertion candidate with zero cost (line~\ref{alg:analyze_candidate:150}).
Otherwise, the final task overlaps with $S.\sigma^{tasks}_h$. 
In this case, in line~\ref{alg:analyze_candidate:160}, we continue to find a feasible insertion candidate by pushing $j$ while keeping its initial task before $S.\sigma^{tasks}_{pos_{t_1}}$.
Finally, the resulting insertion candidate is returned in line~\ref{alg:analyze_candidate:170}.


\begin{algorithm}[!ht]
    \scriptsize
    \caption{\scriptsize ANALYZE\_INSERTION\_CANDIDATE($j$, $pos_{t_1}$, $S$, $I$)}
    \label{alg:analyze_candidate}
    $prev \gets S.\sigma^{tasks}_{pos_{t_1}-1}$\; \label{alg:analyze_candidate:010}
    $st_{t_1} \gets S.\tau_{prev} + p_{prev}$\; \label{alg:analyze_candidate:020}
    \If{$j$'s initial task overlaps  $S.\sigma^{tasks}_{pos_{t_1}}$ when starting at $st_{t_1}$}{ \label{alg:analyze_candidate:030}
        \Return $\mathcal{C}: \langle j, \, \KwFalse, \, pos_{t_1}, \, -, \, st_{t_1}, \,\infty \rangle$\; \label{alg:analyze_candidate:040}
    }
    $st_{t_2} \gets st_{t_1} + a_j + L_j$\; \label{alg:analyze_candidate:050}
    $pos_{t_2} \gets pos_{t_1}$\; \label{alg:analyze_candidate:060}
    \If{$j$'s final task fits before $S.\sigma^{tasks}_{pos_{t_2}}$ when starting at $st_{t_2}$}{ \label{alg:analyze_candidate:070}
        \Return $\mathcal{C}: \langle j, \, \KwTrue, \, pos_{t_1}, \, pos_{t_2}, \, st_{t_1}, \, 0 \rangle$\; \label{alg:analyze_candidate:080}
    }
    \If{$st_{t_2} \geq S.C_{max}$}{ \label{alg:analyze_candidate:090}
        $pos_{t_2} \gets |S.\sigma^{tasks}| + 1$\;\label{alg:analyze_candidate:100}
        $cost \gets st_{t_2} + b_j - S.C_{max}$\; \label{alg:analyze_candidate:110}
        \Return $\mathcal{C}: \langle j, \, \KwTrue, \, pos_{t_1}, \, pos_{t_2}, \, st_{t_1}, \, cost \rangle$\; \label{alg:analyze_candidate:120}
    }
    $h \gets$ Find the expected position for $j$'s final task in the interval [$pos_{t_1}+1$, $|S.\sigma^{tasks}|$]\; \label{alg:analyze_candidate:130}
    \If{$j$'s final task fits before $S.\sigma^{tasks}_{h}$ when starting at $st_{t_2}$}{ \label{alg:analyze_candidate:140}
        \Return $\mathcal{C}: \langle j, \, \KwTrue, \, pos_{t_1}, \, h, \, st_{t_1}, \, 0 \rangle$\;  \label{alg:analyze_candidate:150}
    }
    
    $\mathcal{C} \gets$ FIND\_CANDIDATE\_WITH\_PUSH($j$, $h$, $st_{t_1}$, $st_{t_2}$, $pos_{t_1}$, $S$, $I$)\; \label{alg:analyze_candidate:160}
    \Return $\mathcal{C}$\; \label{alg:analyze_candidate:170}
\end{algorithm}


\begin{algorithm}[!ht]
    \scriptsize
    \caption{\scriptsize FIND\_CANDIDATE\_WITH\_PUSH($j$, $h$, $st_{t_1}$, $st_{t_2}$, $pos_{t_1}$, $S$, $I$)}
    \label{alg:find_cand_w_push}
    $idle \gets$ idle time after $j$'s initial task\; \label{alg:find_cand_w_push:030}
    \While{$h < |S.\sigma^{tasks}|$}{ \label{alg:find_cand_w_push:040}
        $curr \gets S.\sigma^{tasks}_{h}$\; \label{alg:find_cand_w_push:050}
        $push_j \gets S.\tau_{curr} + p_{curr} - st_{t_2}$\; \label{alg:find_cand_w_push:060}
        \If{$push_j > idle$}{ \label{alg:find_cand_w_push:070}
            \Return $\mathcal{C}: \langle j, \, \KwFalse, \, pos_{t_1}, \, h, \, st_{t_1}, \, \infty \rangle$\; \label{alg:find_cand_w_push:080}  
        }
        \If{$j$'s final task fits before $S.\sigma^{tasks}_{h+1}$ when starting at $st_{t_2}+push_j$}{ \label{alg:find_cand_w_push:090}
            \Return $\mathcal{C}: \langle j, \, \KwTrue, \, pos_{t_1}, \, h+1, \, st_{t_1} + push_j, \, 0 \rangle$\;  \label{alg:find_cand_w_push:100}
        }
        $h \gets h + 1$\; \label{alg:find_cand_w_push:150}
    }
    $push_j \gets S.C_{max} - st_{t_2}$\; \label{alg:find_cand_w_push:160}
    \If{$push_j \leq idle$}{ \label{alg:find_cand_w_push:170}
        \Return $\mathcal{C}: \langle j, \, \KwTrue, \, pos_{t_1}, \, h+1, \, st_{t_1} + push_j, \, p_{2j} \rangle$\; \label{alg:find_cand_w_push:180}
    }
    \Return $\mathcal{C}: \langle j, \, \KwFalse, \, pos_{t_1}, \, h, \, st_{t_1}, \, \infty \rangle$\; \label{alg:find_cand_w_push:190}
\end{algorithm}

Algorithm~\ref{alg:find_cand_w_push} presents the pseudocode for FIND\_CANDIDATE\_WITH\_PUSH, called in line~\ref{alg:analyze_candidate:160} of Algorithm~\ref{alg:analyze_candidate}.
Assuming that $j$'s final task overlaps task $S.\sigma^{tasks}_h$ when starting at $st_{t_2}$, it tries to find a feasible insertion candidate for $j$ by pushing its initial and final tasks from their respective start times $st_{t_1}$ and $st_{t_2}$.
Moreover, $j$'s initial task must still fit before $S.\sigma^{tasks}_{pos_{t_1}}$.
In line~\ref{alg:find_cand_w_push:030}, we define a variable $idle$ as the idle time after the initial task of $j$.
It indicates the maximum push allowed on $j$ to fit its initial task at $pos_{t_1}$.
In the while loop of lines~\ref{alg:find_cand_w_push:040}-\ref{alg:find_cand_w_push:150}, we iteratively analyze whether it is possible to skip the overlapping task $S.\sigma^{tasks}_h, \ h < |S.\sigma^{tasks}|$.
In line~\ref{alg:find_cand_w_push:050}, we store in variable $curr$ the current task $S.\sigma^{tasks}_h$.
After that, in line~\ref{alg:find_cand_w_push:060}, we calculate $push_j$, which is the necessary push on $j$ for its final task to skip task $curr$.
If $push_j > idle$ (line~\ref{alg:find_cand_w_push:070}), there is no feasible insertion candidate to insert $j$ at $pos_{t_1}$ and we return an infeasible insertion candidate with infinity cost (line~\ref{alg:find_cand_w_push:080}).
Otherwise, the initial task of job $j$ still fits at $pos_{t_1}$ when starting at $st_{t_1} + push_j$.
If $j$'s final task scheduled at $st_{t_2}+push_j$ fits before $S.\sigma^{tasks}_{h+1}$ (line~\ref{alg:find_cand_w_push:090}), then we return a feasible insertion candidate with zero cost (line~\ref{alg:find_cand_w_push:100}).
In line~\ref{alg:find_cand_w_push:150}, we iterate the variable $h$.
At the end of the while loop (lines~\ref{alg:find_cand_w_push:040}-\ref{alg:find_cand_w_push:150}), $S.\sigma^{tasks}_h$ is the last task of $S.\sigma^{tasks}$.
Therefore, we calculate $push_j$ (line~\ref{alg:find_cand_w_push:160}) and then, if $push_j \leq idle$ (line~\ref{alg:find_cand_w_push:170}), we return a feasible insertion candidate with cost $p_{2j}$ (line~\ref{alg:find_cand_w_push:180}).
Otherwise, we return an infeasible insertion candidate with infinity cost (line~\ref{alg:find_cand_w_push:190}).


\subsubsection{Time complexity}
\label{sec:time_complexity}

We argue that our proposed first-fit algorithm is a polynomial-time heuristic with running time in $O(n^3)$, where $n$ is the number of jobs.
Note that lines~\ref{alg:analyze_candidate:130} and~\ref{alg:analyze_candidate:160} of Algorithm~\ref{alg:analyze_candidate} describe complementary searches, where line~\ref{alg:analyze_candidate:130} loops over $[pos_{t_1}, h)$, and line~\ref{alg:analyze_candidate:160} continues over $[h, |S.\sigma^{tasks}|)$.
For that reason, ANALYZE\_INSERTION\_CANDIDATE has a linear time complexity on the number of jobs, i.e., $O(n)$.
\okC{In lines~\ref{alg:single_dec:035}-\ref{alg:single_dec:100} of Algorithm~\ref{alg:single_dec}, we run ANALYZE\_INSERTION\_CANDIDATE for each position after the finish of $\sigma^{jobs}_{i-1}$ until the end of $\sigma^{tasks}$.
Thus, the time complexity of lines~\ref{alg:single_dec:035}-\ref{alg:single_dec:100} is $O(n^2)$.
Additionally, the time complexity of lines~\ref{alg:single_dec:103}-\ref{alg:single_dec:115} in Algorithm~\ref{alg:single_dec} is $O(n)$, because it applies two linear insertions in the array $S.\sigma^{tasks}$.
As the lines~\ref{alg:single_dec:035}-\ref{alg:single_dec:120} of Algorithm~\ref{alg:single_dec} are executed for each job of $\sigma^{jobs}$, the total time complexity of the \mbox{FIRST\_FIT\_ALGORITHM} is~$O(n^3)$.}

%% file: 173_brkga_single_decoder_first_fit_analysis.tex
\subsubsection{First-fit \okC{quality} analysis \okC{and alternative strategies}}
\label{sec:first_fit_analysis}

Remember that finding an optimal solution \okC{for the SMCTSP} is NP-hard even if the sequence for the initial tasks is predefined \cite{ConSha12}.
Accordingly, given the optimal sequence for the initial tasks, our first-fit approach does not guarantee the generation of an optimal solution.
\okC{We present two instance examples that support the statement, in addition to discussing alternative strategies to generate the optimal schedule for them.}
    
\begin{figure}[!h]
    \centering
    \includegraphics[width=.80\linewidth, trim={1.45cm 0 3cm 0cm}, clip]{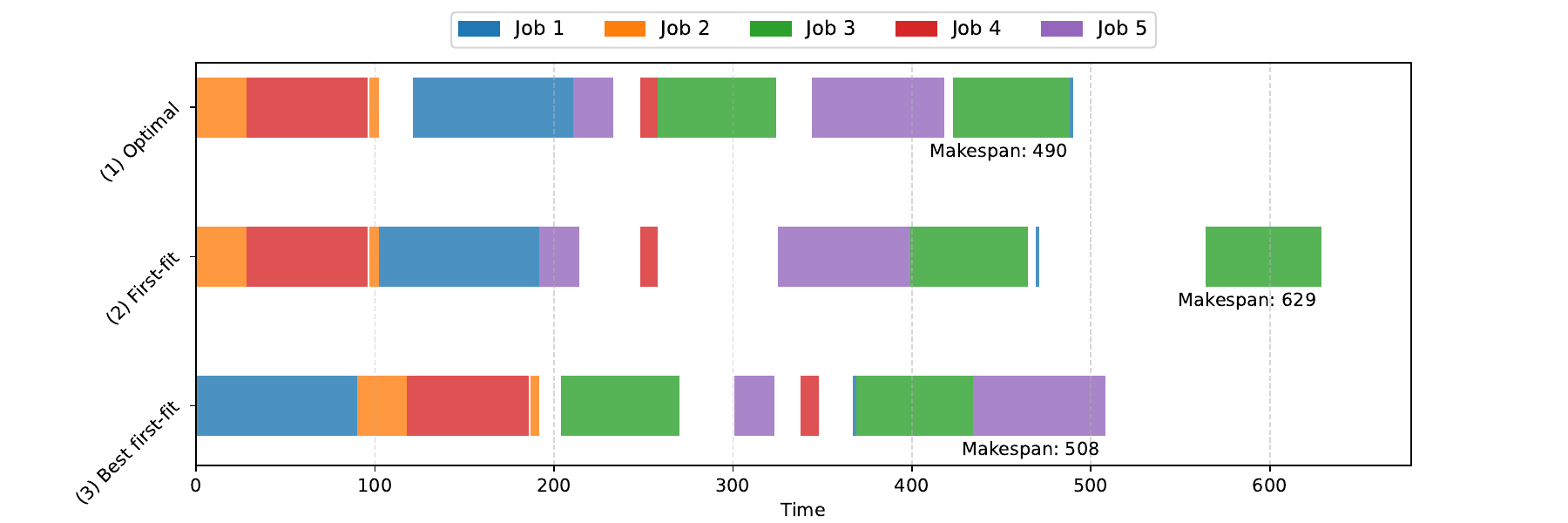}
    \caption{\okC{Schedule comparison for instance 5\_4\_L\_gen of the benchmark set of \protect\citeA{csbenchmark, KhaSalChe20}.}}
    \label{fig:schedule_comparison_ff_vs_opt_scenario_1}
\end{figure}

\okC{Figure \ref{fig:schedule_comparison_ff_vs_opt_scenario_1} exhibits three schedules for the instance \mbox{5\_4\_L\_gen} of the benchmark set of \citeA{csbenchmark, KhaSalChe20}.
\texttt{(1) Optimal} is the optimal schedule, which follows the sequence $\sigma^{jobs} = (2, 4, 1, 5, 3)$ for the initial tasks, and the makespan is equal to 490.
This same sequence $\sigma^{jobs}$ passed to our first-fit algorithm yields a different schedule, referred to as \texttt{(2) First-fit} in Figure \ref{fig:schedule_comparison_ff_vs_opt_scenario_1}, with makespan equal to~629.
Despite this difference being significant, our first-fit algorithm can produce several other schedules.
Given that the instance 5\_4\_L\_gen contains five jobs, there are $5! = 120$ possible sequences of jobs, meaning that our first-fit algorithm can output 120 different schedules.
The best one over them is the schedule \texttt{(3)~Best first-fit} in Figure \ref{fig:schedule_comparison_ff_vs_opt_scenario_1}, which comes from the sequence $\sigma^{jobs}~=~(1,2,4,3,5)$.
Its makespan is equal to $508$, $121$ less than the makespan of the schedule \texttt{(2)~First-fit} in Figure \ref{fig:schedule_comparison_ff_vs_opt_scenario_1}, and only $18$ \okC{($\approx3.7$\%)} more than the makespan of the schedule \texttt{(1) Optimal} in Figure~\ref{fig:schedule_comparison_ff_vs_opt_scenario_1}.}

\okC{Observe that, compared to schedule \texttt{(1)~Optimal} in Figure~\ref{fig:schedule_comparison_ff_vs_opt_scenario_1}, the schedule \texttt{(2)~First-fit} in Figure~\ref{fig:schedule_comparison_ff_vs_opt_scenario_1} has the same task subsequence before the insertion of $\sigma^{jobs}_{5} = 3$, i.e., $\sigma^{tasks} = (3, 7, 4, 1, 9, 8, 10, 2)$ is the task subsequence in both schedules disregarding job~$3$.
Our first-fit algorithm does not schedule job~$3$ before the final task of job~$5$ because its final task would overlap with the final task of job~$1$.
To accommodate job~$3$ at its expected position, another strategy is necessary.}

\begin{figure}[!h]
    \centering
    \includegraphics[width=.80\linewidth, trim={0cm 0 3cm 0cm}, clip]{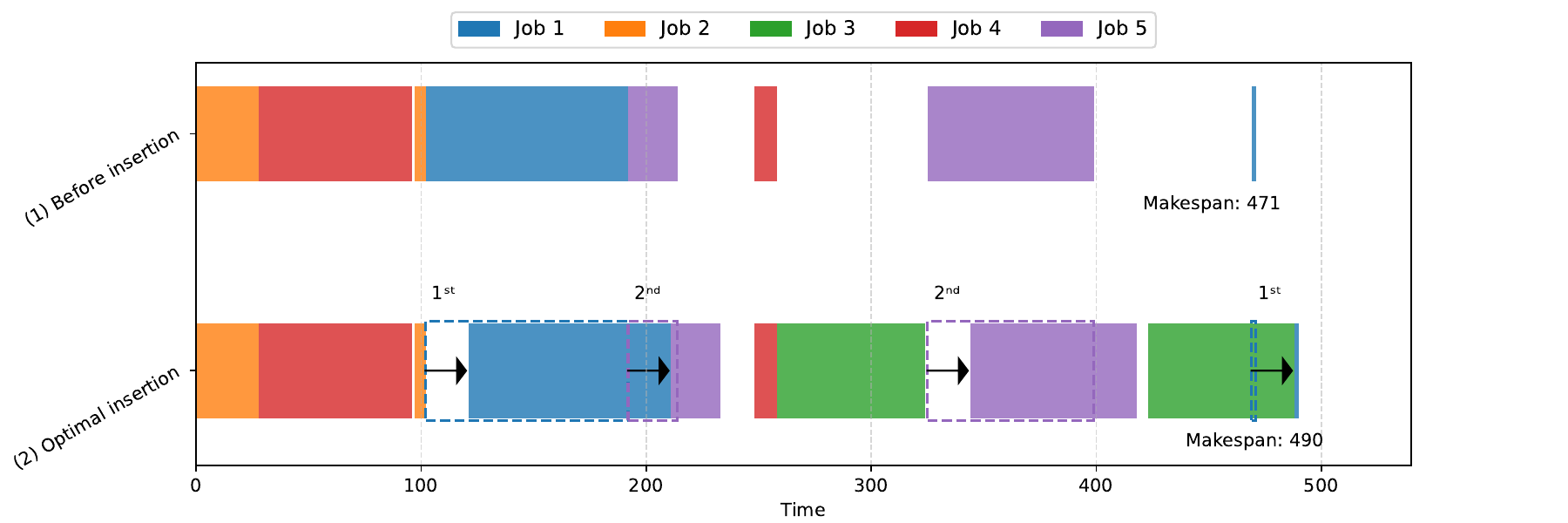}
    \caption{\okC{Optimal insertion for instance 5\_4\_L\_gen of the benchmark set of \protect\citeA{csbenchmark, KhaSalChe20}. Optimality requires right-shifting jobs 1 and 5.}}
    \label{fig:optimal_insertion_strategy_A}
\end{figure}

\okC{Figure~\ref{fig:optimal_insertion_strategy_A} illustrates a possible strategy to obtain the optimal solution for the instance \mbox{5\_4\_L\_gen}.
Shortly, this strategy modifies our first-fit algorithm to allow right-shifting already scheduled jobs in case of overlaps while maintaining the current task sequence.
\texttt{(1) Before insertion} is the schedule achieved by our first-fit algorithm before inserting the job~$3$ in the sequence $\sigma = (2,4,1,5,3)$, which follows the task subsequence $\sigma^{tasks} = (3,7,4,1,9,8,10,2)$.
\texttt{(2)~Optimal insertion} displays the optimal schedule and the modifications applied to the schedule \texttt{(1)~Before insertion} to obtain the optimal insertion for the job~$3$.
An arrow represents the right shift a task received from its previous allocated time, represented by the dashed lines.
The ordinal numbers above the arrows denote the sequence of right-shifts on the respective task.
Note that, instead of skipping the final task of job~$1$, we right-shift it to accommodate job~$3$. 
Consequently, it would also right-shift the job~$5$ because the initial task of job~$1$ would overlap with the initial task of job~$5$.}

\begin{figure}[!ht]
    \centering
    \includegraphics[width=.80\linewidth, trim={1.45cm 0 3cm 0cm},clip]{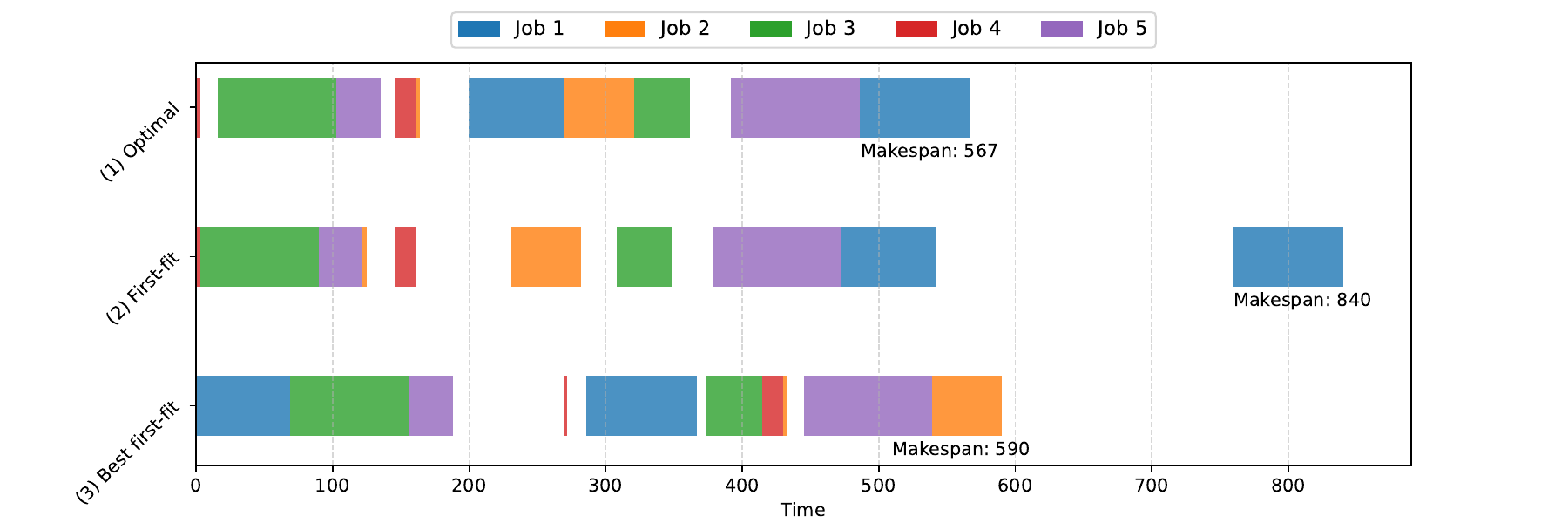}       
    \caption{\okC{Schedule comparison for instance 5\_10\_L\_gen \protect\cite{csbenchmark, KhaSalChe20}.}}
    \label{fig:schedule_comparison_ff_vs_opt_scenario_2}
\end{figure}

\okC{Notice that the strategy illustrated in Figure \ref{fig:optimal_insertion_strategy_A} did not modify the sequence $\sigma^{tasks}$ of the scheduled tasks to accommodate job~$3$.
Even though it is already computationally expensive due to a possible ejection chain, right-shifting scheduled tasks while maintaining the current task sequence to accommodate a job remains a suboptimal approach.
Figure \ref{fig:schedule_comparison_ff_vs_opt_scenario_2} illustrates another instance on which the previous strategy would not produce the optimal solution.
As for Figure~\ref{fig:schedule_comparison_ff_vs_opt_scenario_1}, we present three equivalent schedules, but for the instance 5\_10\_L\_gen of the benchmark set of \citeA{csbenchmark, KhaSalChe20}.}

\okC{In Figure~\ref{fig:schedule_comparison_ff_vs_opt_scenario_2}, the schedule \texttt{(1) Optimal} follows the sequence $\sigma^{jobs} = (4,3,5,2,1)$ for the initial tasks, and the makespan equals $567$.
Similar to Figure \ref{fig:schedule_comparison_ff_vs_opt_scenario_1}, this same sequence $\sigma^{jobs}$ passed to our first-fit algorithm outputs a different schedule, referred to as \texttt{(2) First-fit} in Figure \ref{fig:schedule_comparison_ff_vs_opt_scenario_2}, with makespan equal to $840$.
Nevertheless, similarly to instance 5\_4\_L\_gen, our first-fit algorithm can output $120$ different schedules for instance 5\_10\_L\_gen.
The best one over them is the schedule \texttt{(3)~Best first-fit} in Figure~\ref{fig:schedule_comparison_ff_vs_opt_scenario_2}, which comes from the sequence $\sigma^{jobs} = (1, 3, 5, 4, 2)$.
Its makespan is equal to $590$, $250$ less than the makespan of the schedule \texttt{(2) First-fit}, and only $23$ \okC{($\approx4.1$\%)} more than the makespan of the schedule \texttt{(1)~Optimal}.}

\okC{Notice that, disregarding job~$1$, the task subsequence of the schedule \texttt{(1) Optimal} in Figure~\ref{fig:schedule_comparison_ff_vs_opt_scenario_2} is $\bar{\bar{\sigma}}^{tasks} = (7, 5, 9, 8, 3, 4, 6, 10)$, while the task subsequence of the schedule \texttt{(2) First-fit} in Figure~\ref{fig:schedule_comparison_ff_vs_opt_scenario_2} is $\bar{\sigma}^{tasks} = (7, 5, 9, 3, 8, 4, 6, 10)$.
Remark that the initial task of job~$2$ (task $3$) is placed at position $5$ in schedule \texttt{(1) Optimal}, while this same task is placed at position~$4$ in schedule \texttt{(2) First-fit}.
This suggests another strategy different from the one depicted in Figure~\ref{fig:optimal_insertion_strategy_A}.}



\begin{figure}[!h]
    \centering
    \includegraphics[width=.80\linewidth, trim={0cm 0 3cm 0cm}, clip]{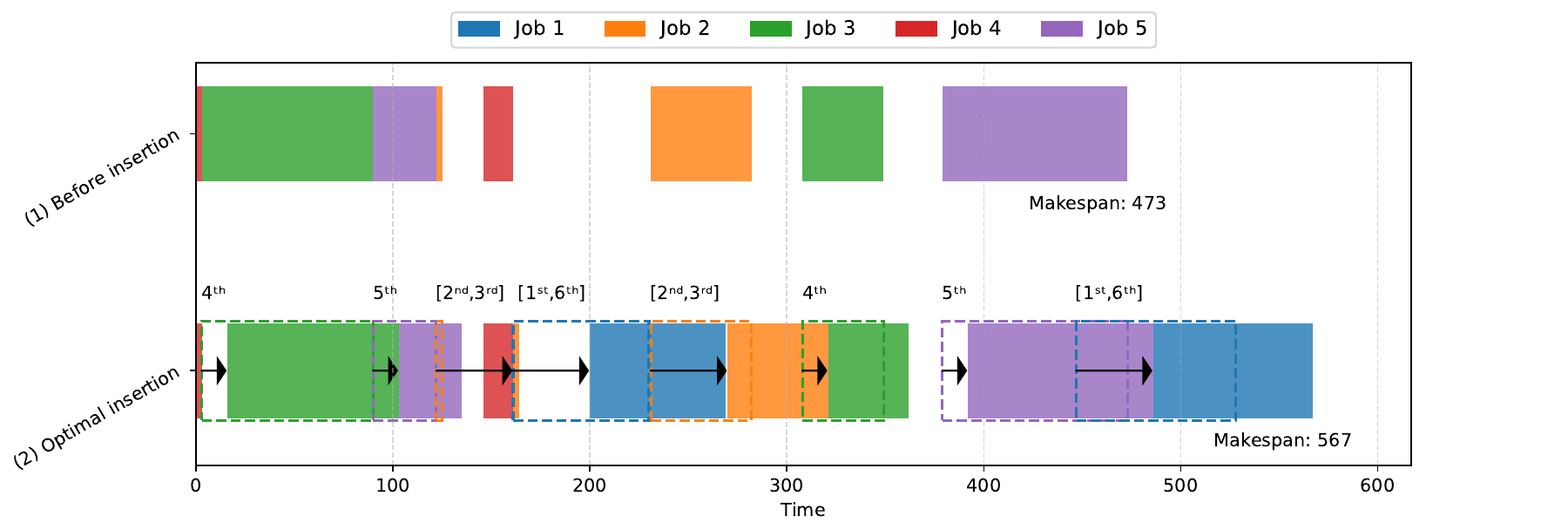}
    \caption{\okC{Optimal insertion for instance 5\_10\_L\_gen of the benchmark set of \protect\citeA{csbenchmark, KhaSalChe20}. Optimality requires right-shifting jobs 1, 2, 3, and 5. The initial task of job 2 skips the final task of job 4.}}
    \label{fig:optimal_insertion_strategy_B}
\end{figure}

\okC{Figure~\ref{fig:optimal_insertion_strategy_B} shows another strategy to obtain the optimal solution for the instance 5\_10\_L\_gen.
Similar to Figure~\ref{fig:optimal_insertion_strategy_A}, the arrows, dashed lines, and the ordinal numbers above the arrows detail the strategy, which allows modifications on the task sequence $\bar{\sigma}^{tasks}$.
In this scenario, the first four jobs of $\sigma^{jobs}$ would be inserted the same way as our first-fit algorithm.
For the fifth job ($\sigma^{jobs}_5 = 1$), however, the new strategy would move the initial task of job~$2$ to position~$5$ after a series of right-shifts.
We don't analyze the step-by-step process to insert job~$1$, although the reader might come up with the same solution based on the arrows, dashed lines, and the ordinal numbers denoting the sequence of the right-shifts.
Jobs with two ordinal numbers mean that they were right-shifted twice.}

\okC{Even though the two presented strategies may yield optimal schedules for a larger set of instances, their computational complexity is significantly greater than our first-fit algorithm.
Furthermore, this trade-off might not be worth it in practice when using them in BRKGA's decoder.
The presented examples show that the best first-fit schedules are reasonably good approximations to the optimal schedules, and the BRKGA embedded with such a decoder will be able to explore the set of possible solutions more efficiently.}

%% file: 175_brkga_single_warm_start.tex
\subsection{\okC{Initial solution generator}}
\label{sec:warm-start}

In a standard BRKGA, the initial population contains only random individuals.
\citeA{LONDE20251} cited several authors who show the benefits of providing good solutions as individuals in the initial population.
In this work, we set part of the initial population with the best solutions of a multi-start procedure.
We first execute an adaptive version of the decoder, described in Algorithm~\ref{alg:first_fit_adap}.
Thereafter, we perform multiple executions of its randomized version, described in Algorithm~\ref{alg:first_fit_adap_rand}.
The best $m \leq p$ solutions are converted into individuals of the initial population.
The remaining $p-m$ individuals are generated randomly as standard.
Preliminary experiments showed that passing a percentage of the multi-start solutions can be beneficial.
Thus, we set $m=\lambda^{ws}\cdot p$, where $\lambda^{ws}$ is an input parameter defining the percentage of the initial population filled with the best multi-start solutions.


Algorithm~\ref{alg:first_fit_adap} describes the pseudocode for the FIRST\_FIT\_ADAPTIVE called in the multi-start procedure.
Lines \ref{alg:first_fit_adap:010}-\ref{alg:first_fit_adap:020} start the partial solution $S$ by inserting the job with the longest delay at the beginning of the schedule.
Lines \ref{alg:first_fit_adap:030}-\ref{alg:first_fit_adap:040} keep track of the last job inserted in $S$ and create a set $R$ of the remaining jobs to be inserted, i.e., not yet in the solution.
After that, the for loop of lines \ref{alg:first_fit_adap:050}-\ref{alg:first_fit_adap:120} schedules the remaining jobs in $R$.
Line \ref{alg:first_fit_adap:060} initializes the best insertion candidate $B\mathcal{C}$ with infinity cost to serve as a sentinel, as this is a minimization problem.
Next, the for loop of lines \ref{alg:first_fit_adap:070}-\ref{alg:first_fit_adap:100} looks for the least cost first-fit insertion candidate among all remaining jobs $j \in R$.
\okC{As in lines \ref{alg:single_dec:040}-\ref{alg:single_dec:100} of Algorithm~\ref{alg:single_dec}, we also execute ANALYZE\_INSERTION\_CANDIDATE($j$, $pos_{t_1}$, $S$, $I$) (Algorithm~\ref{alg:analyze_candidate}) in line \ref{alg:first_fit_adap:080} of Algorithm~\ref{alg:first_fit_adap} for every position $pos_{t_1}$ after the $last_{job}$'s initial task in schedule $S.\sigma^{tasks}$ to insert job $j$.} 
In lines \ref{alg:first_fit_adap:090} and \ref{alg:first_fit_adap:100}, $B\mathcal{C}$ is updated every time $\mathcal{C}.cost < B\mathcal{C}.cost$.
At the end of the for loop of lines \ref{alg:first_fit_adap:070}-\ref{alg:first_fit_adap:100}, 
we update the partial solution $S$ \okC{as in lines \ref{alg:single_dec:103}-\ref{alg:single_dec:115} of Algorithm~\ref{alg:single_dec}, i.e.,} by inserting $B\mathcal{C}.j$ according to  $B\mathcal{C}$ (line \ref{alg:first_fit_adap:105}).
Subsequently, we update $last_{job}$ (line \ref{alg:first_fit_adap:110}) and remove it from the remaining jobs $R$.
Finally, after the for loop of lines \ref{alg:first_fit_adap:050}-\ref{alg:first_fit_adap:120}, $S$ contains a feasible schedule for all jobs $j \in J$ and we return $S$ (line \ref{alg:first_fit_adap_rand:130}).

\begin{algorithm}[H]
    \scriptsize
    \setstretch{0.75}
    \caption{\scriptsize FIRST\_FIT\_ADAPTIVE($I$)}
    \label{alg:first_fit_adap}
    $j \gets$ job with the longest delay\; \label{alg:first_fit_adap:010}
    $S \gets$ Partial solution with the $j$ scheduled at time 0\; \label{alg:first_fit_adap:020}
   $last_{job} \gets j$\; \label{alg:first_fit_adap:030}
    $R \gets$ Copy of $J \setminus \{last_{job}$\}\; \label{alg:first_fit_adap:040}
    \For{$i \gets 2$ \KwTo $n$}{ \label{alg:first_fit_adap:050}
        $B\mathcal{C} \gets \mathcal{C} : \langle 0, \KwFalse, 0, 0, 0, \infty \rangle$\; \label{alg:first_fit_adap:060}
        \For{$j \in R$}{\label{alg:first_fit_adap:070}
            $\mathcal{C} \gets$ First-fit insertion candidate after $last_{job}$'s initial task in $S.\sigma^{tasks}$ to insert $j$\; \label{alg:first_fit_adap:080}
            \If{$\mathcal{C}.cost < B\mathcal{C}.cost$}{ \label{alg:first_fit_adap:090}
                $B\mathcal{C} \gets \mathcal{C}$\;\label{alg:first_fit_adap:100}
            }
        }
        $S \gets$ Update partial solution $S$ by inserting $B\mathcal{C}.j$ according to $B\mathcal{C}$\; \label{alg:first_fit_adap:105}
        $last_{job} \gets B\mathcal{C}.j$\; \label{alg:first_fit_adap:110}
        Remove $last_{job}$ from $R$\; \label{alg:first_fit_adap:120}
    }
    
    \Return $S$\; \label{alg:first_fit_adap:130}
\end{algorithm}


Algorithm~\ref{alg:first_fit_adap_rand} describes the randomized version of the FIRST\_FIT\_ADAPTIVE presented in Algorithm~\ref{alg:first_fit_adap}.
The main differences are highlighted in the following.
First, in line \ref{alg:first_fit_adap_rand:010}, the first job of the schedule is selected randomly.
Then, we initialize an empty insertion candidate list (line \ref{alg:first_fit_adap_rand:060}).
It will store all the first-fit insertion candidates $\mathcal{C}$ of the jobs $j \in R$ (line \ref{alg:first_fit_adap_rand:090}).
Afterward, an insertion candidate is chosen from $CL$ by randomly selecting an insertion candidate $\mathcal{C}$ from $CL$ restricted by quality, which requires a threshold parameter $\alpha \in [0, 1]$.
For this, in line \ref{alg:first_fit_adap_rand:095}, we create a restricted insertion candidate list $RCL$ composed of insertion candidates $\mathcal{C} \in CL$, whose cost $\mathcal{C}.cost$ lies in the interval $[c_{min}, c_{min} + \alpha (c_{max} - c_{min})]$, where $c_{min} = \texttt{min}(\mathcal{C}.cost : \mathcal{C} \in CL)$, and  $c_{max} = \texttt{max}(\mathcal{C}.cost : \mathcal{C} \in CL)$.
Later, in line \ref{alg:first_fit_adap_rand:100}, we randomly choose an insertion candidate $\mathcal{C}$ from $RCL$.
Note that for $\alpha = 0$, the FIRST\_FIT\_ADAPTIVE\_RANDOMIZED is not equivalent to the FIRST\_FIT\_ADAPTIVE of Algorithm~\ref{alg:first_fit_adap} because two insertion candidates in $RCL$ can have the same cost, i.e., the algorithm can randomly choose between them, instead of selecting the first one found.

\begin{algorithm}[!ht]
    \scriptsize
    \setstretch{0.75}
    \caption{\scriptsize FIRST\_FIT\_ADAPTIVE\_RANDOMIZED($I$, $\alpha$)}
    \label{alg:first_fit_adap_rand}
    $j \gets$ random job\; \label{alg:first_fit_adap_rand:010}
    $S \gets$ Partial solution with the $j$ scheduled at time 0\; \label{alg:first_fit_adap_rand:020}
    $last_{job} \gets j$\; \label{alg:first_fit_adap_rand:030}
    $R \gets$ Copy of $J \setminus \{last_{job}\}$\; \label{alg:first_fit_adap_rand:040}
    \For{$i \gets 2$ \KwTo $n$}{\label{alg:first_fit_adap_rand:050}
        $CL \gets \varnothing$\; \label{alg:first_fit_adap_rand:060}
        \For{$j \in R$}{\label{alg:first_fit_adap_rand:070}
            $\mathcal{C} \gets$ First-fit insertion candidate after $last_{job}$'s initial task in $S.\sigma^{tasks}$ to insert $j$\; \label{alg:first_fit_adap_rand:080}
            Add $\mathcal{C}$ to $CL$\; \label{alg:first_fit_adap_rand:090}
        }
        $RCL \gets$ $CL$ filtered by quality using threshold parameter $\alpha$\; \label{alg:first_fit_adap_rand:095}
        $\mathcal{C} \gets$ Random candidate from $RCL$\; \label{alg:first_fit_adap_rand:100}
        $S \gets$ Update partial solution $S$ by inserting $\mathcal{C}.j$ according to $\mathcal{C}$\; \label{alg:first_fit_adap_rand:110}
        $last_{job} \gets \mathcal{C}.j$\; \label{alg:first_fit_adap_rand:120}
        Remove $last_{job}$ from $R$\; \label{alg:first_fit_adap_rand:130}
    }
    
    \Return $S$\; \label{alg:first_fit_adap:140}
\end{algorithm}

%% file: 180_brkga_single_local_search.tex
\subsection{Local search}
\label{sec:local_search}

\citeA{LONDE20251} presented an extensive list of authors that applied local search (LS) algorithms in the BRKGA.
They hybridize the LS with the BRKGA in several ways.
In this work, we combine two of these ideas.
First, inspired by \citeA{ANDRADE201967}, we periodically apply the local search to the best individual of the elite set.
Second, inspired by \citeA{silva2024}, we apply the local search after an improvement on the best $b$ individuals of the elite set.
Besides, our LS method is a modification of the commonly used \textit{insert} method \cite{den2001neighborhoods}.
Our modification allows more control over the neighborhood size.
In this paper, we adopted the term \textit{move} instead of \textit{insert} for this method.
In what follows, we describe our neighborhood structure design (Section \ref{sec:ls_neighborhood}), present the local search algorithm (Section \ref{sec:ls_algorithm}), and detail the integration with the BRKGA (Section \ref{sec:ls_integration_brkga}).

%% file: 182_brkga_ls_neighborhood.tex
\subsubsection{Local search neighborhood}
\label{sec:ls_neighborhood}

Remember that in Section \ref{sec:encoding_decoder_brkga} we stated that our first-fit algorithm is injective.
Thus, we represent each solution by the sequence $\sigma^{jobs}$ retrieved from a vector of random keys $X$.
Given a sorted sequence of the jobs $\sigma^{jobs}$ of an individual of the current BRKGA population, define the \textit{move} operation as one that removes the job in position $i$ and inserts it in position $j$ of $\sigma^{jobs}$.
A neighbor of $\sigma^{jobs}$ is a sorted sequence that can be obtained from $\sigma^{jobs}$ through a \textit{move} operation.
The size of this neighborhood is known to be $(n-1)^2$ \cite{den2001neighborhoods}.

\begin{figure}[h]
    \centering
    \includegraphics[width=0.90\linewidth]{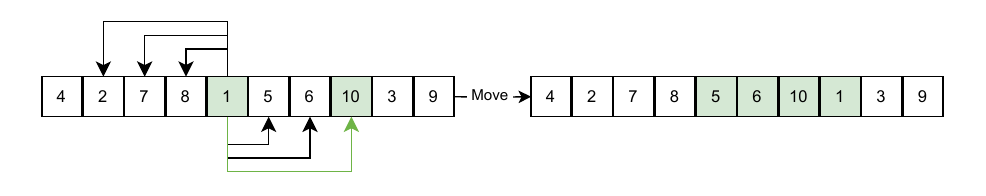}
    \caption{Move candidates for $\sigma^{jobs}_5 = 1$ with radius $r=3$ and a move operation.}
    \label{fig:move_candidates_job}
\end{figure}

Notice that the makespan calculus of each neighbor could lead to an ejection chain that reconstructs the entire solution.
This ejection chain is mainly caused by interleaving jobs, i.e., when tasks of different jobs are scheduled between a job's tasks.
Therefore, we run the FIRST\_FIT\_ALGORITHM described in Algorithm \ref{alg:single_dec}, which has a time complexity of $O(n^3)$ (see Section \ref{sec:time_complexity}).
Consequently, the time complexity of one LS iteration is $O(n^5)$.
Given that we execute the LS periodically, this method becomes computationally expensive.
On that account, we propose a parameter that controls the neighborhood size.
Given a sorted sequence $\sigma^{jobs}$, define the radius $r$ as the extension of the \textit{move} operation, i.e., $\sigma^{jobs}_k$ ($1 \leq k \leq n$) can only be moved to position $l$ (max$(1, k-r) \leq l \leq$ min$(k+r, \ n)$).
Figure \ref{fig:move_candidates_job} illustrates the \textit{move} candidates of a job within a radius and the result of a \textit{move} operation.

Given that we apply the LS in two different situations, we denote by $r_{pLS}$ and $r_{iLS}$ the radius values for the periodic LS and the LS after an improvement, respectively.
Note that, if the radius $r$ of the LS iteration is bounded above by a small constant, the number of neighbors becomes $O(n)$, which also decreases the time complexity of one LS iteration from $O(n^5)$ to $O(n^4)$.

%% file: 183_brkga_ls_algorithm.tex
\subsubsection{Local search algorithm}
\label{sec:ls_algorithm}

Algorithm \ref{alg:move_job_first_improv} describes the local search algorithm, which receives as input a random-key vector $X$, its makespan $C_{max}$, the local search radius $r$ and the instance $I$.
In lines \ref{alg:move_job_first_improv:010}-\ref{alg:move_job_first_improv:020}, we initialize the best random-key vector and its makespan.
Next, we first set a flag $impr$ as \KwTrue (line \ref{alg:move_job_first_improv:030}) to enter the while loop (line \ref{alg:move_job_first_improv:040}), then we set it to \KwFalse (line \ref{alg:move_job_first_improv:050}).
If an improving neighbor is found during the while loop of lines \ref{alg:move_job_first_improv:040}-\ref{alg:move_job_first_improv:220}, $impr$ becomes \KwTrue and we initiate a new LS iteration (line \ref{alg:move_job_first_improv:040}).
Otherwise, we leave the while loop because we reached a local optimum.
In line \ref{alg:move_job_first_improv:060}, we generate the sorted sequence $\sigma^{jobs}$ of the random-key vector $X^*$ (see Section \ref{sec:encoding_decoder_brkga}).
In line \ref{alg:move_job_first_improv:070}, we set a new variable $k$ as $1$, referencing the job $\sigma^{jobs}_k$ to be moved.

Afterward, in the while loop of the lines \ref{alg:move_job_first_improv:080}-\ref{alg:move_job_first_improv:220}, we analyze, for each $k \in \{1,\ldots,n\}$, each move candidate of the job $\sigma^{jobs}_k$ (see Figure \ref{fig:move_candidates_job}) until an improvement is found (i.e., $impr =$ \KwTrue) or all jobs have been analyzed (i.e., $k > n$).
In line \ref{alg:move_job_first_improv:090}, we initialize variable $l$ as max($1$, $k-r$), which stores the destination candidate of the job $\sigma^{jobs}_k$.
In line \ref{alg:move_job_first_improv:100}, we move the job $\sigma^{jobs}_k$ to position $l$ in $\sigma^{jobs}$ to start the analysis of the move candidates. 

In the while loop of the lines \ref{alg:move_job_first_improv:110}-\ref{alg:move_job_first_improv:180}, we increment over the extension of the move operation (i.e., $l \leq k+r$) until we reach the end of the schedule (i.e., $l \leq n$) or find an improvement (i.e., $impr =$ \KwTrue).
In line \ref{alg:move_job_first_improv:120}, we check whether the candidate move destination is valid by ensuring it is neither the origin (i.e., $l \neq k$) nor a previously analyzed neighbor (i.e., $l \neq k - 1$). 
Consequently, when this condition is met, line \ref{alg:move_job_first_improv:130} decodes $\sigma^{jobs}$ to yield the corresponding solution $S$.
If the makespan of $S$ is smaller than the makespan of $X^*$ (line \ref{alg:move_job_first_improv:140}), an improvement was found.
Therefore, in line \ref{alg:move_job_first_improv:150}, we update the best makespan $C^*_{max}$, the best random-key vector $X^*$, and activate the flag $impr$.
Then, if the next job destination is within the schedule (line \ref{alg:move_job_first_improv:160}), we swap $\sigma^{jobs}_l$ with $\sigma^{jobs}_{l+1}$ (line \ref{alg:move_job_first_improv:170}), i.e., the moving job that was in position $l$ is now placed in position $l+1$.
Next, we increment $l$ (line \ref{alg:move_job_first_improv:180}).

At the end of the while loop of lines \ref{alg:move_job_first_improv:110}-\ref{alg:move_job_first_improv:180}, there are two possible outcomes.
First, all move candidates were analyzed, and no improvement was found.
In this case (line \ref{alg:move_job_first_improv:190}), we return the moving job that is in position min($l$, $n$) (line \ref{alg:move_job_first_improv:200}) to its original position $k$ in $\sigma^{jobs}$ (line \ref{alg:move_job_first_improv:210}) and we increment $k$ to analyze the next move candidates for job $\sigma^{jobs}_{k+1}$ (line \ref{alg:move_job_first_improv:220}).
Otherwise, an improvement was found and $\sigma^{jobs}$ is the improving neighbor.
Therefore, we leave the loop of lines \ref{alg:move_job_first_improv:110}-\ref{alg:move_job_first_improv:220} and start a new local search round.
After the loop of lines \ref{alg:move_job_first_improv:040}-\ref{alg:move_job_first_improv:220},
$X^*$ encodes a local optimum solution and $C^*_{max}$ corresponds to its makespan.
Thus, in line \ref{alg:move_job_first_improv:230}, we return both variables.

\begin{algorithm}[!ht]
    \scriptsize
    \setstretch{0.75}
    \caption{\scriptsize MOVE\_JOB\_FIRST\_IMPROVEMENT($X$, $C_{max}$, $r$, $I$)}
    \label{alg:move_job_first_improv}
    $X^* \gets$ copy of $X$\; \label{alg:move_job_first_improv:010}
    $C^*_{max} \gets C_{max}$\; \label{alg:move_job_first_improv:020}
    $impr \gets \KwTrue$\; \label{alg:move_job_first_improv:030}
    \While{$imprv$}{ \label{alg:move_job_first_improv:040}
        $impr \gets \KwFalse$\; \label{alg:move_job_first_improv:050}
        $\sigma^{jobs} \gets$ sorted sequence of the jobs of the vector $X^*$\; \label{alg:move_job_first_improv:060}
        $k \gets 1$\; \label{alg:move_job_first_improv:070}
        \While{$k \leq n$ \KwAnd \KwNot $imprv$}{ \label{alg:move_job_first_improv:080}
            $l \gets$ max($1$, $k-r$)\; \label{alg:move_job_first_improv:090}  
            $\sigma^{jobs} \gets$ Move $\sigma^{jobs}_k$ to position $l$ in $\sigma^{jobs}$\; \label{alg:move_job_first_improv:100}
            \While{$l \leq$ min($k+r, \ n$) \KwAnd \KwNot $impr$}{\label{alg:move_job_first_improv:110}
                \If{$l \neq k$ \KwAnd $l \neq k - 1$}{ \label{alg:move_job_first_improv:120}
                    $S \gets$ FIRST\_FIT\_ALGORITHM($\sigma^{jobs}$, $I$)\; \label{alg:move_job_first_improv:130}
                    \If{$S.C_{max} < C^*_{max}$}{ \label{alg:move_job_first_improv:140}
                        Update $C^*_{max}$, $X^*$ and set $impr \gets \KwTrue$\; \label{alg:move_job_first_improv:150}
                    }
                }
                \If{$l + 1 \leq n$}{ \label{alg:move_job_first_improv:160}
                    $\sigma^{jobs} \gets$ Swap $\sigma^{jobs}_l$ with $\sigma^{jobs}_{l+1}$\; \label{alg:move_job_first_improv:170}
                }
                $l \gets l + 1$\; \label{alg:move_job_first_improv:180}
            }
            \If{\KwNot $impr$}{ \label{alg:move_job_first_improv:190}
                $l \gets$ min($l$, $n$)\; \label{alg:move_job_first_improv:200}
                $\sigma^{jobs} \gets$ Move $\sigma^{jobs}_l$ to position $k$ in $\sigma^{jobs}$\;     \label{alg:move_job_first_improv:210}
                $k \gets k + 1$\; \label{alg:move_job_first_improv:220}
            }
        }
    }
    \Return $X^*$, $C^*_{max}$ \label{alg:move_job_first_improv:230}
\end{algorithm}

%% file: 184_brkga_ls_integration.tex
\subsubsection{Integration with the BRKGA}
\label{sec:ls_integration_brkga}

We propose to execute the LS in two situations: (i) periodically and (ii) after an improvement.
The periodic LS is set to run every $L$ iterations on the best \okC{eligible} individual of the current elite set.
After each periodic LS, we mark the individual as ineligible for another periodic LS.
We also set $L$ to be proportional to the number of jobs $n$.
The LS after an improvement occurs on the best $b$ \okC{eligible} individuals of the elite set.
Similar to the periodic LS, we mark each individual as ineligible for any other LS.
Note that, for $r_{pLS}<r_{iLS}$, the LS after an improvement can still be executed on the individual that has already undergone a periodic LS. 
In turn, the periodic LS is not applied \okC{to the individual on whom the LS was performed} after an improvement.
Besides, given a population $P$ sorted by fitness, whenever we search for the next best individual, we store the fitness of the previously analyzed individual $i \in P$.
If the next individual $i+1 \in P$ has the same fitness, we skip the LS execution. 

\okC{Table \ref{tab:eligibility_ls} illustrates an example of LS eligibility within the elite set. 
Column ``Pos.'' shows each individual’s position.
Column ``Eligibility'' lists four cases: 
``Ineligible (I)'' when LS was applied after an improvement; 
``Ineligible (P)'' when a periodic LS was applied; 
``Eligible'' when no LS was applied;
and ``Same fitness'' when its fitness matches that of someone ranked above.
Finally, column ``Fitness'' exhibits the fitness of the individual. Assuming $b=2$, observe that the LS after an improvement would be executed on the second and third individuals of the elite set, while the periodic LS would be executed only on the fourth individual.}

\begin{table}[!htp]
\centering
\caption{\okC{Example of LS eligibility within the elite set}}
\label{tab:eligibility_ls}
\scriptsize
\setlength{\tabcolsep}{8pt}
\renewcommand{\arraystretch}{1.3}

\definecolor{mygreen}{HTML}{b6d7a8}
\definecolor{myorange}{HTML}{f9cb9c}
\definecolor{myblue}{HTML}{a4c2f4}
\definecolor{myred}{HTML}{ea9999}
\definecolor{myyellow}{HTML}{ffe599}

\begin{tabular}{|l|l|r|}
\hline
\textbf{Pos.} & \textbf{Eligibility} & \textbf{Fitness} \\
\hline
\rowcolor{myred}     1  & Ineligible (I)      & 4347 \\
\rowcolor{myorange}     2  & Ineligible (P)    & 4524 \\
\rowcolor{myorange}     3  & Ineligible (P)    & 4554 \\
\rowcolor{mygreen}    4  & Eligible     & 4562 \\
\rowcolor{myorange}     5  & Ineligible (P)    & 4599 \\
\rowcolor{mygreen}   6  & Eligible          & 4845 \\
\rowcolor{myyellow}  7  & Same fitness         & 4845 \\
\rowcolor{mygreen}   8  & Eligible          & 4846 \\
\rowcolor{mygreen}   9  & Eligible          & 4860 \\
\rowcolor{mygreen}   10 & Eligible          & 4862 \\
\hline
\end{tabular}
\end{table}

%% file: 185_brkga_single_perturbation.tex
\subsection{Perturbation}
\label{sec:perturbation}

This section describes the perturbation component of our BRKGA.
Throughout the generations, the BRKGA aims to improve the population by combining pairs of solutions from the elite and non-elite sets, but biased towards the elite solutions.
After a certain number of iterations, the BRKGA expects the elite set to contain a set of blocks of keys in the random-key vectors that positively impact the quality of a solution.
We denote this set as the \textit{convergence structure}.
Notwithstanding this, the BRKGA may eventually reach many generations without improving the overall best solution.
The metaheuristics generally contain a diversification mechanism to overcome this situation, which can provide new search directions.
For the BRKGA, it consists of the insertion of randomly generated individuals.
However, sometimes this may not be sufficient to produce improved solutions, and adding a perturbation component such as a restart can be beneficial.

During some preliminary tests, we noticed that parameter configurations can highly affect the number of unique random-key vectors in a population throughout the generations.
In Appendix~\ref{app:pop_homogeneity_brkga}, we show the correlation between the percentage of unique random-key vectors in a population and 
(i) the elite percentage $p_e$, 
(ii) the mutant percentage $p_m$, 
(iii) and the inheritance probability of an elite parent key $\rho_e$.
In summary, BRKGA-R can quickly converge to a homogeneous elite set after each restart. 
Additionally, for higher $p_e$, lower $p_m$, and higher $\rho_e$, part of the non-elite set becomes less diverse.
This is a common concern in genetic algorithms, denoted as \textit{premature convergence} \cite{Fogel1994}, and there are some approaches to prevent it in genetic algorithms \cite{PANDEY20141047}.
\citeA{Londe2024} states that in a BRKGA, the premature convergence is related to a lack of diversity in the elite set, although we also observed it in the non-elite set.
They also review some operators to prevent the premature convergence, including the island model \cite{PANDEY20141047, Whitley1998}, the reset operator, and the shake \cite{ANDRADE201967}.

Given that we observed this issue even with the restarts, a possible workaround would be to reduce the restart periodicity.
The problem with this approach is that each time a population is restarted, it destroys the convergence structure of the population, i.e., valuable parts of the random-key vector that usually appear in better solutions can be lost.
On that account, we adopted the \textit{shaking} method, a BRKGA feature introduced by \citeA{ANDRADE201967} to deal with this scenario.
It randomly modifies the individuals in the elite set and restarts only the non-elite set, instead of restarting the whole population.
Hence, the elite set will preserve some parts of the convergence structure.

\subsubsection{Shaking}
\label{sec:shaking_description}

The original BRKGA shaking method proposed by \citeA{ANDRADE201967} is executed over the problem solution space, i.e., over the decoded solutions.
Some authors followed this approach \cite{LONDE2021114728, Kim01122024}.
On the other hand, other authors have proposed to shake over the BRKGA space, i.e., over the random-key vectors~\cite{LONDE2022109634, mauri2021}.
In this work, we employ the shaking method implemented in the framework of \citeA{ANDRADE202117}, which performs the operations over the BRKGA space.
In summary, the shaking method performs $\psi = \lambda^{shake} \cdot n$ pairs of random operations in each elite set individual and replaces the non-elite set members with random samples.
The parameter $\lambda^{shake} \in [0, 1]$ is a percentage of the number of jobs.
There are two types of shaking ($s^{type}$): \textit{CHANGE} and \textit{SWAP}.
In \textit{CHANGE}, it first selects a random key and inverts its value.
Then, it selects another random key and replaces its value with another random value.
In \textit{SWAP}, it first selects a random key and swaps its value with the following key's value.
In sequence, it selects two random keys and swaps their values.

\begin{figure}[h]
    \centering
    \includegraphics[width=0.48\linewidth, trim={3cm 3.5cm 3.2cm 3.5cm}, clip]{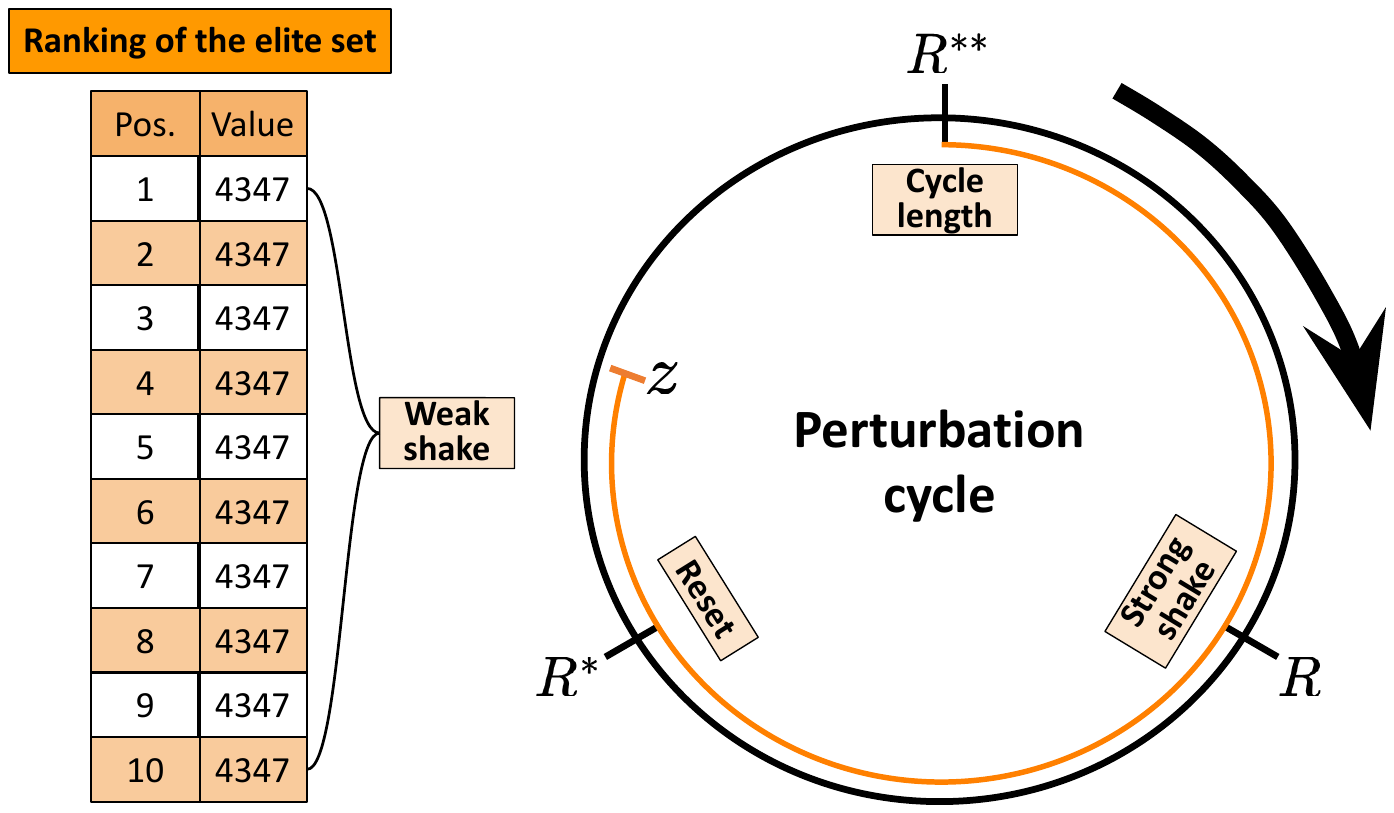}
    \caption{\okC{Perturbation cycle of our BRKGA}}
    \label{fig:perturbation_cycle}
\end{figure}

As the \textit{shaking} method enables more frequent perturbations, it can be applied in more situations.
To address the issue of premature convergence, we apply a weak shake each time the elite set becomes homogeneous, i.e., the best fitness is equal to the worst fitness of the elite set.
We define the weak shake as \citeA{ANDRADE201967}, i.e., $\lambda^{shake}$ is uniformly drawn from the interval $[0.05, 0.2]$.
Sometimes, the BRKGA gets stuck, and the weak shake may not help the BRKGA to find new improved solutions.
In this case, we first define $z$ as the current number of iterations without improvement.
Then, we apply a strong shake each time \okC{$z \mod R^{**} = R$, where $R^{**}$} defines the length of the perturbation cycle and $R$ defines the strong shake iteration in this perturbation cycle.
We define the strong shake as \citeA{ANDRADE201967}, i.e., $\lambda^{shake}$ is uniformly drawn from the interval $[0.5, 1.0]$.
Finally, we still reset the population each time \okC{$z \mod R^{**} = R^*, \ R^* > R$}, where $R^*$ defines the restart iteration in the perturbation cycle.
\okC{Figure \ref{fig:perturbation_cycle} illustrates the perturbation cycle of our BRKGA.}

However, some observations indicated that the BRKGA may lose its convergence structure even when we apply a weak shake.
This is similar to what happens when resetting the population.
Some authors solve this issue by injecting a random-key vector, such as the overall best solution found so far or the best \okC{initial} solution \cite{Festa2024, DEABREU2022101149, ANDRADE201967}, but this has not been applied after the \textit{shaking} to our knowledge.
Nonetheless, preliminary tests showed that injecting a random-key vector can be beneficial.
Accordingly, we selected some injection candidates from the preliminary tests.
The first candidate is the type ``best solution of the current population'', denoted as CB.
It is a good candidate to inject after a weak shake, as we want to preserve the convergence structure.
The second candidate is the type ``best overall solution found so far'', denoted as OB.
It performs well after each type of perturbation, whether it is a weak shake, a strong shake, or a population reset.
The third candidate is the type ``best solution of a new multi-start execution'', denoted as BMS.
Instead of injecting the \okC{initial} solution again, we provide a new random solution with potentially new and useful parts.
We define parameters that indicate the type of solution to be injected after a weak shake ($\gamma^{weak}$), a strong shake ($\gamma^{strong}$), and a reset ($\gamma^{reset}$), where $\gamma^{weak} \in \{$CB, OB$\}$, $\gamma^{strong} \in \{$OB, BMS$\}$, and $\gamma^{reset} \in \{$OB, BMS$\}$.
\okC{Figure~\ref{fig:injection_strategy} displays the injection strategy after perturbations in our BRKGA.}

\input{186_extra_visual_diagram}

%% file: 186_extra_visual_diagram.tex
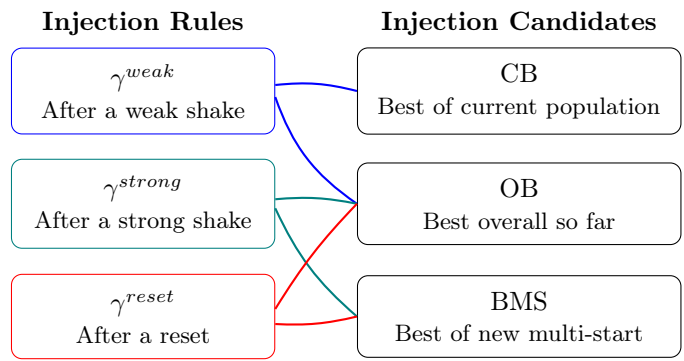
\begin{figure}[ht!]
\centering
\begin{tikzpicture}[
    font=\small,
    box/.style={rectangle, draw, rounded corners, align=center, inner sep=6pt},
    arrow/.style={-, thick, >=latex},
    node distance=1cm
]

\node[align=center] at (5,0.8) {\textbf{Injection Rules}};

\node[box, minimum width=3.5cm, draw=blue] (weak)  at (5,-0.1)
    {$\gamma^{weak}$ \\ \footnotesize After a weak shake};

\node[box, minimum width=3.5cm, draw=teal] (strong) at (5,-1.6)
    {$\gamma^{strong}$ \\ \footnotesize After a strong shake};

\node[box, minimum width=3.5cm, draw=red] (resetp) at (5,-3.1)
    {$\gamma^{reset}$ \\ \footnotesize After a reset};

\node[align=center] at (10,0.8) {\textbf{Injection Candidates}};

\node[box, minimum width=4.3cm] (cb)  at (10,-0.1) 
    {CB \\ \footnotesize Best of current population};

\node[box, minimum width=4.3cm] (ob)  at (10,-1.6) 
    {OB \\ \footnotesize Best overall so far};

\node[box, minimum width=4.3cm] (bms) at (10,-3.1) 
    {BMS \\ \footnotesize Best of new multi-start};


\draw[arrow, blue]
  (weak.east) ++(0, 0.08) to[bend left=10] (cb.west);

\draw[arrow, blue]
  (weak.east) ++(0,-0.08) to[bend right=18] (ob.west);

\draw[arrow, teal]
  (strong.east) ++(0, 0.06) to[bend left=8] (ob.west);

\draw[arrow, teal]
  (strong.east) ++(0,-0.06) to[bend right=12] (bms.west);

\draw[arrow, red]
  (resetp.east) ++(0, 0.10) to[bend left=6] (ob.west);

\draw[arrow, red]
  (resetp.east) ++(0,-0.10) to[bend right=8] (bms.west);

\end{tikzpicture}
\caption{\okC{Injection strategy after perturbations in our BRKGA.}}
\label{fig:injection_strategy}
\end{figure}

%% file: 200_brkga_single_overall_framework.tex
\subsection{Overall framework}
\label{sec:overall_framework}

\begin{figure}[!ht]
    \centering
    \includegraphics[width=.90\linewidth]{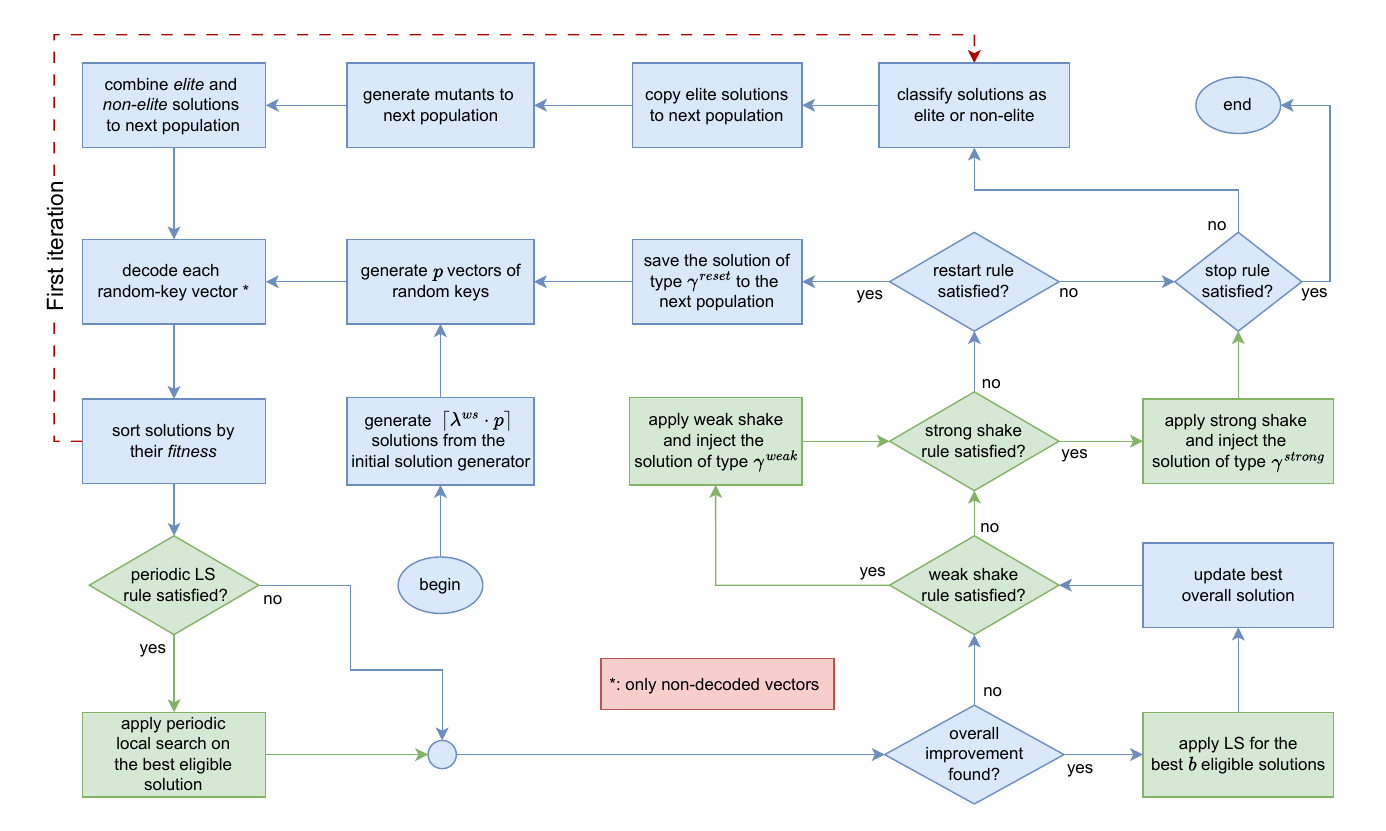}
    \caption{Flowchart of the BRKGA approaches used in this paper for the single-machine CTSP. In blue, the baseline BRKGA-R. In green, the BRKGA-R-S-LS.}
    \label{fig:brkga-flowchart}
\end{figure}

Figure \ref{fig:brkga-flowchart} depicts our proposed BRKGA's flowchart for the single-machine CTSP.
The baseline BRKGA-R is presented in blue, while the additional components of the \mbox{BRKGA-R-S-LS} are highlighted in green.
In the beginning, we generate $\ceil{\lambda^{ws}\cdot p}$ solutions with the \okC{initial solution generator} (see Section \ref{sec:warm-start}), where $p$ is the population size, and $\lambda^{ws}$ is the percentage of these solutions that will be inserted into the initial population.
Then, we encode these solutions and generate other $\floor{(1-\lambda^{ws})\cdot p}$ random-key vectors to construct the initial population.
The random-key vectors of the initial population are decoded and sorted by their fitness.
Afterward, the main loop of the BRKGA starts.
At the beginning of each BRKGA iteration, we compute the next generation of the population based on solution quality.
Consequently, we classify the solutions as elite or non-elite, where the first $\floor{p_e\cdot p}$ solutions compose the elite set.
Then, the next population consists of: the current elite set, a set of $\floor{p_m\cdot p}$ mutant solutions, and the remaining solutions resulting from the combination of elite and non-elite solutions.
Afterward, the non-decoded random-key vectors are decoded and sorted by their fitness.
Given an iteration $i$ of the main BRKGA loop, we execute the periodic LS on the best eligible individual of the population whenever \okC{$i \mod L = 0$} (see Section~\ref{sec:ls_integration_brkga}).
If an overall improvement was found in the current generation, we apply the LS on the best $b$ eligible individuals and update the best overall individual (see Section \ref{sec:ls_integration_brkga}).
After that, we consider three perturbations (see Section \ref{sec:perturbation}).
First, if the condition for a weak shake is satisfied, we apply it and inject the solution of type $\gamma^{weak}$.
Second, if the condition for a strong shake is satisfied, we apply it and inject the solution of type $\gamma^{strong}$.
Third, if the condition for a restart is satisfied, we restart the BRKGA, resetting the population by replacing it with $p-1$ new random-key vectors and the solution of type $\gamma^{reset}$.
While the stopping rule is not satisfied, we repeat the process.
Otherwise, the algorithm finishes, and the best random-key vector obtained is returned.

\vspace{-0.5cm}

%% file: 250_computationalresults.tex
\section{Computational experiments}
\label{sec:experiments}

In this section, we present the computational results to evaluate the performance of the \okC{CP} model and the BRKGA metaheuristic approaches. 
All experiments were conducted on a machine running Ubuntu x86-64 GNU/Linux, with an Intel Core i7-10700 Octa-Core 2.90 GHz processor and 16 GB of RAM.
The \okC{CP} model was implemented in MiniZinc version 2.8.7 \cite{Nethercote2007, Stuckey2014} and solved with OR Tools CP-SAT version 9.11.4210 \cite{cpsatlp}. 
\okC{It is worth highlighting that we use the MiniZinc's predicate \texttt{disjunctive} for the noOverlap global constraint \eqref{cp:03} in the CP model (see Section \ref{sec:notationscp}).}
The BRKGA was implemented in Julia version 1.11.3, based on the BrkgaMpIpr.jl framework~\cite{brkgajulia, ANDRADE202117}.
\okC{We also implemented the state-of-the-art MIP formulations of \citeA{KhaSal24} in Julia 1.11.3, using the following packages: JuMP v1.23.6, Gurobi v1.5.0, Gurobi\_jll v12.0.2.}

\okC{We organize this section as follows.
Section~\ref{sec:instancias} presents the benchmark instances.
Section~\ref{sec:tested_approaches_parameter_settings} lists the tested approaches and details the parameter settings.
Section~\ref{sec:comparing_cp_vs_mips} compares the CP with the state-of-the-art MIP formulations.
Section~\ref{sec:comparison_brkgas_cp_bks} shows a comparison between our approaches and the best-known solutions.
Section~\ref{sec:comparison_tested_approaches} compares the tested approaches with each other under the same computational settings.
Section~\ref{sec:impact_brkga_components} analyzes the impact of the BRKGA components.}

%% file: 255_benchmark_instances.tex
\subsection{Benchmark instances}
\label{sec:instancias}

We considered the general benchmark set proposed in \citeA{csbenchmark, KhaSalChe20}. 
There are eight different numbers of jobs, i.e., $n \in \{5, 10, 15, 20, 25, 40, 50, 100\}$. 
For each value of $n$, there are $30$ randomly generated instances divided into three categories: small (S), medium (M), and large (L).
In total, there are 240 instances. 
Each category has 10 instances, where the processing time of the tasks and the delay durations are drawn from the discrete uniform distribution with the following parameters: for category S, $a_j, b_j \sim U(1, 20)$ and $L_j \sim U(10, 80)$; for category M, $a_j, b_j \sim U(1, 50)$ and $L_j \sim U(25, 200)$; and for category L, $a_j, b_j \sim U(1, 100)$ and $L_j \sim U(50, 400)$. 
The instances are represented as $n$\_$id$\_category\_gen, where $id \in [1,10]$ and “gen” indicates that the instances are part of the general set proposed by \citeA{KhaSalChe20}.

%% file: 260_parametrization.tex
\subsection{Tested approaches and parameter settings}
\label{sec:tested_approaches_parameter_settings}

The following approaches were evaluated in the computational experiments: 
\begin{itemize}
\setlength\itemsep{-.3em}
    \item the \okC{CP} model described in Section \ref{sec:notationscp}, denoted as CP;
    \item the BRKGA with the restart presented in Section \ref{sec:brkga}, denoted as BRKGA-R;
    \item the BRKGA with the restart, the shake, and the local search detailed in Section \ref{sec:brkga}, denoted as BRKGA-R-S-LS;
    \item \okC{the MIP formulation of \citeA{BekGalJunOswRei14}, denoted as MIP\_P1 \cite{KhaSal24};}
    \item \okC{the MIP formulation of \citeA{KhaSal24}, denoted as MIP\_P2.}
\end{itemize}

We execute ten independent runs of each BRKGA metaheuristic in a single thread for three minutes.
\okC{We remark that we don't compare the BRKGAs with the state-of-the-art heuristics of \citeA{KhaSal24} (BS\_P1, BS\_P2, and R\&S) since they are based on mathematical formulations and require significantly more time (60 minutes) to obtain competitive results compared to our BRKGAs, which are designed to run within low computational times (three minutes in our tests).}
We run the OR Tools CP-SAT with two different setups.
The first executes one run in a single thread with a time limit of three minutes, only to compare with the BRKGA metaheuristics.
The second uses eight threads with a time limit of \okC{one hour} \okC{to: (1) find optimal solutions; (2) compute the best objective values and lower bounds; and (3) compare with the state-of-the-art MIP formulations}.
Both setups stop when they find an optimal solution, or they reach their respective time limits.
In terms of notation, we additionally add a suffix to differentiate the setups.
CP-1T denotes the CP executed with one thread and a time limit of three minutes, while CP-8T denotes the CP executed with eight threads and a time limit of one hour.
\okC{For the state-of-the-art formulations MIP\_P1 and MIP\_P2, we execute one run of the solver using eight threads with a time limit of one hour to compare with the CP-8T.}

\okC{We have many parameters to set up for the two BRKGA approaches we implemented (see Section~\ref{sec:brkga}).}
Therefore, we employed the \texttt{irace} package version $4.2.0$ developed by \citeA{iracepackage} to define the parameters for each algorithm.
Table \ref{tab:summary_all_parameters} presents the parameter settings suggested by \texttt{irace}.
We assume $r_{iLS} = n$, i.e., we search the whole neighborhood for the LS after an improvement.
Besides, we also assume $\gamma^{reset} \in \{$OB, \okC{BI}, BMS$\}$ for \okC{the BRKGA-R}, where \okC{BI} is the type ``best \okC{initial} solution''.
It is worth mentioning that setting $\gamma^{reset} =$ CB is equivalent to assuming $\gamma^{reset} = $ OB for the BRKGA-R.
For details on the parameters and configurations, see Appendix~\ref{app:parameter_settings}.

\begin{table}[!htp]\centering
\caption{\okC{Summary of the parameters, their definitions, and the settings suggested by \texttt{irace}}}
\label{tab:summary_all_parameters}
\scriptsize
\setstretch{0.9}
\okC{
\resizebox{0.95\textwidth}{!}{
\begin{tabular}{l m{1cm} p{6cm} l c c l c c}
\toprule
\multirow{2}{*}{Component} &
\multirow{2}{*}{Notation} &
\multirow{2}{*}{Definition} & &
\multicolumn{2}{c}{BRKGA-R} & &
\multicolumn{2}{c}{BRKGA-R-S-LS} \\
\cmidrule{5-6} \cmidrule{8-9}
& & & & Used & \texttt{irace} & & Used & \texttt{irace} \\
  \cmidrule{1-3} \cmidrule{5-6} \cmidrule{8-9}  
\multirow{5}{*}{Standard}
& $n$ & number of random keys in a solution & & \checkmark & --- & & \checkmark & --- \\
& $p$ & population size & & \checkmark & 163 & & \checkmark & 185 \\
& $p_e$ & elite population size & & \checkmark & 0.39 & & \checkmark & 0.43 \\
& $p_m$ & mutant population size & & \checkmark & 0.20 & & \checkmark & 0.24 \\
& $\rho_e$ & probability of inheritance of an elite parent key & & \checkmark & 0.74 & & \checkmark & 0.78 \\
  \cmidrule{1-3} \cmidrule{5-6} \cmidrule{8-9}  
\multirow{3}{*}{$\parbox{2cm}{Initial solution\\generator}$}
& $\alpha$ & threshold parameter for quality selection & & \checkmark & 0.02 & & \checkmark & 0.01 \\
& $\lambda^{ws}$ & percentage of the initial population from multi-start & & \checkmark & 0.28 & & \checkmark & 0.22 \\
& $n_{msi}$ & number of multi-start iterations & & \checkmark & 1374 & & \checkmark & 602 \\
 \cmidrule{1-3} \cmidrule{5-6} \cmidrule{8-9} 
\multirow{8}{*}{Perturbation}
& $n_{nimp}$ & restart cycle while no improvement & & \checkmark & 956 & --- & --- \\
& $R$ & strong shake iteration in the perturbation cycle & & --- & --- & & \checkmark & 154 \\
& $R^{*}$ & restart iteration in the perturbation cycle & & --- & --- & & \checkmark & $2\cdot R$ \\
& $R^{**}$ & length of the perturbation cycle & & --- & --- & & \checkmark & $9\cdot R$ \\
& $s^{type}$ & shake type — \{CHANGE, SWAP\} & & --- & --- & & \checkmark & SWAP \\
& $\gamma^{weak}$ & type of solution injected after a weak shake & & --- & --- & & \checkmark & OB \\
& $\gamma^{strong}$ & type of solution injected after a strong shake & & --- & --- & & \checkmark & OB \\
& $\gamma^{reset}$ & type of solution injected after a reset & & \checkmark & BMS & & \checkmark & OB \\
  \cmidrule{1-3} \cmidrule{5-6} \cmidrule{8-9}  
\multirow{4}{*}{Local search}
& $b$ & max.\ number of elite-set individuals for LS \linebreak after improvement & & --- & --- & & \checkmark & 9 \\
& $L$ & LS periodicity & & --- & --- & & \checkmark & $0.21\cdot n$ \\
& $r_{pLS}$ & LS radius for periodic LS & & --- & --- & & \checkmark & 7 \\
& $r_{iLS}$ & LS radius for LS after improvement & & --- & --- & & \checkmark & --- \\
\bottomrule
\end{tabular}
}}
\end{table}

%% file: 263_comparison_cp_vs_mips_local.tex
\subsection{\okC{Comparison between the CP and the state-of-the-art MIP formulations}}
\label{sec:comparing_cp_vs_mips}

\okC{In this section, we compare the CP-8T against the MIP\_P1 and MIP\_P2 of \citeA{KhaSal24} within the same computational environment (see Section \ref{sec:tested_approaches_parameter_settings}).
The MIP approaches obtained feasible solutions for 230 (MIP\_P1) and 206 (MIP\_P2) out of the 240 instances (95.83\% and 85.83\%).
\okC{We highlight that MIP\_P2 did not obtain a feasible solution for any instance with 100 jobs.}
Nonetheless, our CP-8T obtained feasible solutions for all instances.
Besides, the CP-8T and the MIP approaches achieved optimal solutions for every instance with up to ten jobs in at most 30 seconds.
Since no further improvements are possible for these instances with up to ten jobs, we disregard their results.
Additionally, we omit the execution time of the three approaches because almost all of them reached 3600 seconds for the considered instances (except some instances with 15~jobs).}

\okC{Table \ref{tab:comparison_cp_vs_mips} presents an overview of the CP-8T results compared to the MIP\_P2 and MIP\_P1 for $n \geq 15$, grouped by their number of jobs $n$ and their category \texttt{cat}.
We denote each group as \mbox{\texttt{<n> x <cat>}}, where \texttt{<n>} is the number of jobs and \texttt{<cat>} is the category of the instances.
We consider the following metrics over each group: the mean objective value, the mean optimality gap (in percentage), and the number of optimal solutions. 
The optimality gap is calculated as: $100 \cdot\frac{|value - bound|}{value}$, where $bound$ is the objective bound, and $value$ is the objective value.
The value ``NA'' (Not Applicable) is placed whenever the solver on the respective formulation did not find a feasible solution for at least one instance of the group.
The row \texttt{Total} summarizes the results of all the instances by summing the columns.
The row \texttt{Mean} summarizes the results of all the instances by taking the mean of the columns\okC{, but only for the applicable rows without ``NA'' in the three formulations.}}

\begin{table}[!ht]
\scriptsize
\centering
\caption{\okC{Overview of the CP-8T results compared to the MIP\_P1 and MIP\_P2}}
\label{tab:comparison_cp_vs_mips}
\okC{
\resizebox{\textwidth}{!}{
\begin{tabular}{rl l rrr l rrr l ccc}
\toprule
\multirow{2}{*}{n} & \multirow{2}{*}{cat} &  & \multicolumn{3}{c}{Objective value} & & \multicolumn{3}{c}{Optimality gap (\%)} & & \multicolumn{3}{c}{Number of optimal solutions} \\ 
\cmidrule{4-6} \cmidrule{8-10} \cmidrule{12-14} 
&  &  & CP-8T & MIP\_P2 & MIP\_P1 &  & CP-8T & MIP\_P2 & MIP\_P1 &  & CP-8T & MIP\_P2 & MIP\_P1 \\
\cmidrule{1-2} \cmidrule{4-6} \cmidrule{8-10} \cmidrule{12-14} 
15 & S &  & \textbf{302.20} & 303.50 & 303.70 &  & \textbf{0.22} & 14.75 & 20.44 &  & \textbf{8} & 4 & 2 \\
15 & M &  & \textbf{785.10} & 787.50 & 789.60 &  & \textbf{1.77} & 20.75 & 26.20 &  & \textbf{1} & \textbf{1} & 0 \\
15 & L &  & \textbf{1543.20} & 1550.20 & 1551.80 &  & \textbf{1.62} & 17.14 & 22.55 &  & 2 & \textbf{3} & \textbf{3} \\
\addlinespace
20 & S &  & \textbf{416.90} & 423.00 & 426.20 &  & \textbf{2.79} & 51.88 & 58.09 &  & 0 & 0 & 0 \\
20 & M &  & \textbf{1052.30} & 1063.90 & 1071.70 &  & \textbf{3.46} & 50.37 & 56.13 &  & 0 & 0 & 0 \\
20 & L &  & \textbf{2121.70} & 2147.20 & 2157.30 &  & \textbf{3.55} & 53.58 & 59.01 &  & 0 & 0 & 0 \\
\addlinespace
25 & S &  & \textbf{541.80} & 554.80 & 559.20 &  & \textbf{2.91} & 68.23 & 71.78 &  & 0 & 0 & 0 \\
25 & M &  & \textbf{1338.10} & 1360.10 & 1376.10 &  & \textbf{3.90} & 66.86 & 71.67 &  & 0 & 0 & 0 \\
25 & L &  & \textbf{2595.50} & 2639.70 & 2670.50 &  & \textbf{3.90} & 67.69 & 71.29 &  & 0 & 0 & 0 \\
\addlinespace
40 & S &  & \textbf{850.10} & 885.10 & 906.70 &  & \textbf{3.21} & 81.41 & 84.97 &  & 0 & 0 & 0 \\
40 & M &  & \textbf{2101.60} & 2167.80 & 2246.30 &  & \textbf{4.37} & 80.62 & 84.71 &  & 0 & 0 & 0 \\
40 & L &  & \textbf{4229.40} & 4371.20 & 4475.90 &  & \textbf{4.40} & 81.21 & 84.80 &  & 0 & 0 & 0 \\
\addlinespace
50 & S &  & \textbf{1071.90} & 1117.70 & 1267.60 &  & \textbf{3.61} & 84.79 & 89.61 &  & 0 & 0 & 0 \\
50 & M &  & \textbf{2608.80} & NA & 3059.00 &  & \textbf{4.30} & NA & 88.83 &  & 0 & 0 & 0 \\
50 & L &  & \textbf{5353.70} & NA & 8373.60 &  & \textbf{4.81} & NA & 90.78 &  & 0 & 0 & 0 \\
\addlinespace
100 & S &  & \textbf{2233.00} & NA & NA &  & \textbf{4.25} & NA & NA &  & 0 & 0 & 0 \\
100 & M &  & \textbf{5366.40} & NA & NA &  & \textbf{5.20} & NA & NA &  & 0 & 0 & 0 \\
100 & L &  & \textbf{10748.70} & NA & NA &  & \textbf{5.51} & NA & NA &  & 0 & 0 & 0 \\ 
\cmidrule{1-2} \cmidrule{4-6} \cmidrule{8-10} \cmidrule{12-14} 
\multicolumn{2}{l}{Total} &  & NA & NA & NA &  & NA & NA & NA &  & \textbf{11} & 8 & 5 \\
\multicolumn{2}{l}{Mean} &  & \textbf{1457.68} & 1490.13 & 1523.28 &  & \textbf{3.05} & 56.87 & 61.64 &  & NA & NA & NA \\ 
\bottomrule
\end{tabular}}}
\end{table}

\okC{We notice from Table \ref{tab:comparison_cp_vs_mips} that the CP-8T obtained the best metrics presented for all groups.
Although the solver on the \mbox{MIP\_P2} struggled to find feasible solutions for larger instances, the \mbox{MIP\_P2} demonstrated superiority over the \mbox{MIP\_P1}.
Furthermore, the \mbox{MIP\_P2} outperformed the CP-8T in objective value in only two instances (each with 15 jobs), while the MIP\_P1 did so in only one instance (also with 15 jobs).
However, these differences are actually minimal in practice.
Moreover, note that the number of optimal solutions in the group \texttt{15 x L} is three for the MIPs, and two for the CP-8T.
Despite these few outliers, we conclude that the CP-8T is superior to the MIP formulations.}

\FloatBarrier

%% file: 272_comparison_our_approaches_vs_bks.tex
\subsection{\okC{Comparison between our approaches and the best-known solutions}}
\label{sec:comparison_brkgas_cp_bks}

\okC{This section analyzes the performance of our main BRKGA approach (BRKGA-R-S-LS) and the CP-8T compared to the current best-known solutions (BKS).
We consider the BKS as the best solution obtained by the MIP formulations and heuristics of \citeA{KhaSal24} (MIP\_P1, \mbox{MIP\_P2}, BS\_P1, BS\_P2, R\&S), and our main approaches (CP-8T, BRKGA-R, BRKGA-R-S-LS).}
\okC{It is worth highlighting that we consider both the MIP solutions reported in the literature and those obtained in this work.}

\begin{figure}[h]
    \centering
    \includegraphics[width=0.45\linewidth]{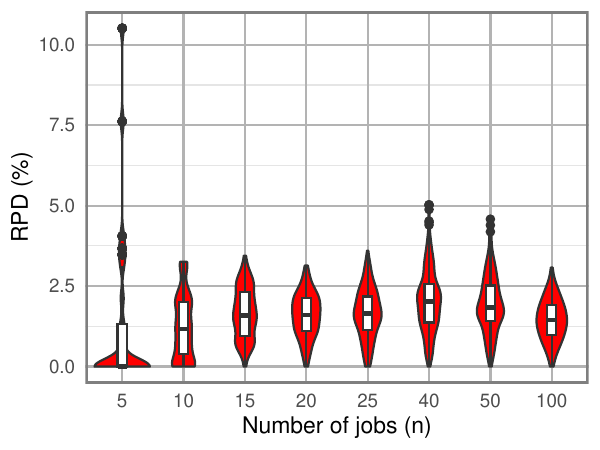}
    \caption{Distribution of RPD values for the BRKGA-R-S-LS compared to the BKS}
    \label{fig:distribution_rpd_brkga_r_s_ls_vs_bks}
\end{figure}

Figure~\ref{fig:distribution_rpd_brkga_r_s_ls_vs_bks} exhibits the distribution of RPD values for the BRKGA-R-S-LS compared to the BKS.
\okC{The Relative Percentage Deviation (RPD) is calculated as $100\cdot \frac{obj^{run} - obj^{\text{BKS}}}{obj^{\text{BKS}}}$, where $obj^{run}$ is the makespan obtained on a run of the BRKGA-R-S-LS, and $obj^{\text{BKS}}$ is the makespan of the BKS.}
Additionally, each box plot contains ten observations for each of the $30$ instances with the same number of jobs $n$, totaling $300$ observations per box plot.
We observe that at least approximately 75\% of the RPD values are at most 2.5\%.
Besides, apart from the groups with five, 40, and 50 jobs, no other group had outliers, and there are no significant differences between the medians, although it indicates a reduction as the number of jobs increases from $n\geq 40$.
Therefore, these results indicate that the BRKGA-R-S-LS is a stable metaheuristic for this problem, obtaining high-quality approximate solutions within low computational times.

\FloatBarrier

\okC{Table~\ref{tab:comparison_cp_bks} presents an overview of the CP-8T results in comparison with BKS for $n \geq 15$ (see Section~\ref{sec:comparing_cp_vs_mips}), grouped by their number of jobs $n$ and their category \texttt{cat}.
Column \texttt{=BKS} shows the number of instances in which CP-8T obtained the BKS value.
Column \texttt{opt} displays the number of proven optimal solutions found by CP-8T.
Column \texttt{RPD (\%)} exhibits the mean RPD value from the BKS, calculated as: $100\cdot\frac{obj^{\text{CP-8T}}-obj^{\text{BKS}}}{obj^{\text{BKS}}}$, where $obj^{\text{CP-8T}}$ is the objective value of the CP-8T, and $obj^{\text{BKS}}$ is the objective value of the BKS.
Column \texttt{gap (\%)} shows the mean gap value of the solver, i.e., $gap = 100 \cdot\frac{|bound - value|}{value}$, where $bound$ is the objective bound, and $value$ is the objective value.
The row \texttt{Total} summarizes the results of all the instances by summing the columns.
The row \texttt{Mean} summarizes the results of all the instances by taking the mean of the columns.}

\begin{table}[!ht]
\scriptsize
\setstretch{1.0}
\centering
\caption{\okC{Overview of the CP-8T results in comparison with the BKS}\label{tab:comparison_cp_bks}}
\okC{\begin{tabular}{ll l r r r r}
  \toprule
  n & cat & & =BKS & opt & RPD $\pm$ std (\%) & gap $\pm$ std (\%) \\ 
  \cmidrule{1-2} \cmidrule{4-7}
  \multirow{3}{*}{15} & S & & 10 & 8 & 0.00 $\pm$ 0.00 & 0.22 $\pm$ 0.59 \\ 
   & M & & 9 & 1 & 0.06 $\pm$ 0.18 & 1.77 $\pm$ 1.07 \\ 
   & L & & 9 & 2 & 0.01 $\pm$ 0.02 & 1.62 $\pm$ 1.28 \\ 
  \addlinespace
  \multirow{3}{*}{20} & S & & 7 & 0 & 0.14 $\pm$ 0.24 & 2.78 $\pm$ 0.49 \\ 
   & M & & 10 & 0 & 0.00 $\pm$ 0.00 & 3.46 $\pm$ 0.16 \\ 
   & L & & 10 & 0 & 0.00 $\pm$ 0.00 & 3.55 $\pm$ 0.34 \\ 
  \addlinespace
  \multirow{3}{*}{25} & S & & 9 & 0 & 0.02 $\pm$ 0.05 & 2.91 $\pm$ 0.54 \\ 
   & M & & 9 & 0 & 0.03 $\pm$ 0.09 & 3.89 $\pm$ 0.28 \\ 
   & L & & 8 & 0 & 0.07 $\pm$ 0.18 & 3.90 $\pm$ 0.44 \\ 
  \addlinespace
  \multirow{3}{*}{40} & S & & 10 & 0 & 0.00 $\pm$ 0.00 & 3.21 $\pm$ 0.49 \\ 
   & M & & 7 & 0 & 0.14 $\pm$ 0.28 & 4.37 $\pm$ 0.90 \\ 
   & L & & 9 & 0 & 0.03 $\pm$ 0.09 & 4.40 $\pm$ 0.61 \\ 
  \addlinespace
  \multirow{3}{*}{50} & S & & 10 & 0 & 0.00 $\pm$ 0.00 & 3.61 $\pm$ 0.40 \\ 
   & M & & 8 & 0 & 0.07 $\pm$ 0.15 & 4.30 $\pm$ 0.46 \\ 
   & L & & 10 & 0 & 0.00 $\pm$ 0.00 & 4.81 $\pm$ 0.42 \\ 
  \addlinespace
  \multirow{3}{*}{100} & S & & 10 & 0 & 0.00 $\pm$ 0.00 & 4.25 $\pm$ 0.34 \\ 
   & M & & 9 & 0 & 0.05 $\pm$ 0.17 & 5.21 $\pm$ 0.42 \\ 
   & L & & 9 & 0 & 0.04 $\pm$ 0.14 & 5.51 $\pm$ 0.39 \\ 
  \cmidrule{1-2} \cmidrule{4-7}
  \multicolumn{2}{l}{Total} & & 163 & 11 & --- & --- \\
  \multicolumn{2}{l}{Mean} & & --- & --- & 0.04 $\pm$ 0.13 & 3.54 $\pm$ 1.41 \\ 
  \bottomrule
\end{tabular}}
\end{table}

\okC{We notice from Table \ref{tab:comparison_cp_bks} that the CP-8T reached the BKS value for 163 out of the 180 instances (90.56\%).}
Regarding optimality, the CP-8T obtained optimal solutions for 11 out of the 180 instances (6.11\%), all containing 15 jobs. 
The RPD value is only 0.04\% on average with a standard deviation of 0.13\%, i.e., despite the CP-8T not obtaining the BKS in some instances, the solution values obtained are similar.
Finally, the solver gap is only 3.54\% on average with a standard deviation of 1.41\%, reaching a maximum average gap of 5.51\% with a standard deviation of 0.39\% for the group \mbox{\texttt{100 x L}}.

\FloatBarrier

%% file: 265_comparison_tested_approaches.tex
\subsection{\okC{Comparison between our proposed approaches}}
\label{sec:comparison_tested_approaches}

This section evaluates the proposed approaches by comparing their performance regarding solution quality \okC{given the same computational environment}.
\okC{Hence}, we only consider the following approaches: BRKGA-R, BRKGA-R-S-LS, and CP-1T (see Section~\ref{sec:tested_approaches_parameter_settings}).
Moreover, this analysis encompasses all the 240 instances of the benchmark, given that the three selected approaches did not obtain the optimal solution for all the instances with up to ten jobs.

\begin{figure}[!ht]
    \centering
    \subfigure[All instances\label{fig:violin_plot_all_insts}]{
        \includegraphics[width=.4\linewidth]{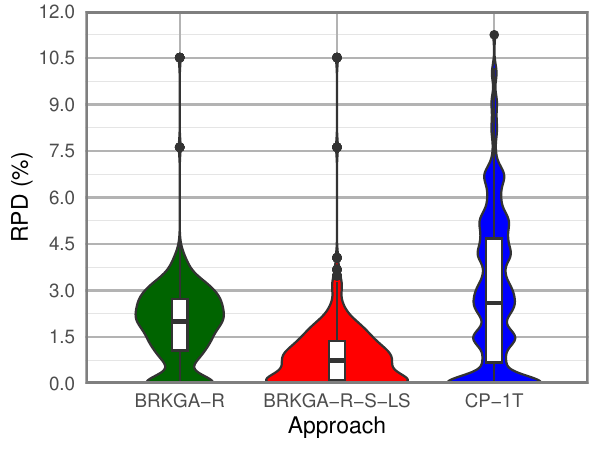}
    }
    \subfigure[Instances with $n\geq 15$ \label{fig:violin_plot_insts_n_min_15}]{
        \includegraphics[width=.4\linewidth]{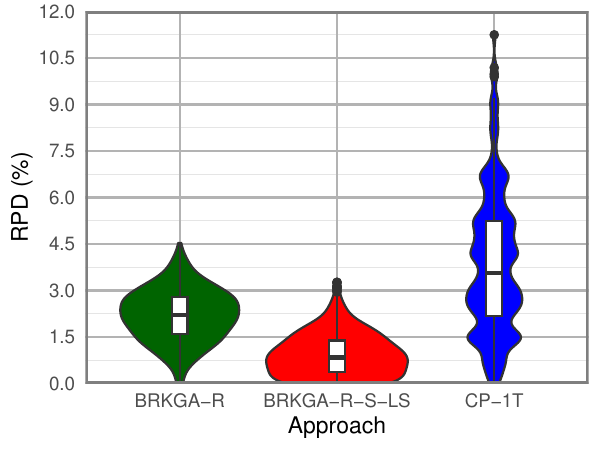}
    }
    \caption{Distribution of RPD values of each approach.}
    \label{fig:distribution_rpd_algorithms}
\end{figure}

Figure \ref{fig:distribution_rpd_algorithms} depicts a distribution of the \okC{RPD} values for two subsets of the benchmark, where each violin plot has a bandwidth equal to $0.25$.
\okC{The RPD is calculated as $100 \cdot \frac{obj^{run} - obj^{best}}{obj^{best}}$, where $obj^{run}$ is the makespan obtained on a run of an approach and $obj^{best}$ is the best makespan found across the three approaches.}
The graphics in Figure \ref{fig:violin_plot_all_insts} consider the distributions for all 240 benchmark instances.
The categories BRKGA-R and BRKGA-R-S-LS  contain one observation for each of the ten independent runs, totaling 2400 observations.
Conversely, category CP-1T contains 240 observations, given that the CP-1T was executed only once for each instance.
Comparison of the medians indicates the superior performance of BRKGA-R-S-LS.
Furthermore, CP-1T indicated the highest elongated distribution, suggesting inconsistent performance.
In contrast, BRKGA-R-S-LS had the narrowest distribution, establishing it as the most stable and robust method.
The bump arising in the lower tail of the distributions indicates that although BRKGA-R-S-LS obtained the highest number of solutions equal to the best solution among the three approaches, many of these solutions were also obtained by the other approaches.
Besides, we highlight the presence of outliers above 7.5\% in categories BRKGA-R and BRKGA-R-S-LS.

Figure \ref{fig:violin_plot_insts_n_min_15} filters the distribution of RPD values of the three approaches for the 180 instances with at least 15 jobs.
Hence, there are 1800 observations for BRKGA-R and BRKGA-R-S-LS, and 180 observations for CP-1T.
This filtered distribution emphasizes the superior performance of \mbox{BRKGA-R-S-LS} compared to BRKGA-R, and in particular to CP-1T.
Notice that for at least 75\% of the observations of BRKGA-R-S-LS, the RPDs were smaller than at least 75\% of the observations of the other approaches.
Furthermore, for all observations of BRKGA-R-S-LS, the RPDs were smaller than at least 50\% of the observations of the CP-1T.
For at least 75\% of the observations of BRKGA-R, the RPDs were smaller than at least 50\% of the RPDs of the CP-1T.
Moreover, a direct comparison between Figures \ref{fig:violin_plot_all_insts} and~\ref{fig:violin_plot_insts_n_min_15} shows that the bump in the lower tail of each category and the outliers in the BRKGAs disappeared in the filtered distribution.
These facts suggest that many of the observations that obtained solutions equal to the best solution among the three approaches and the outliers of the BRKGAs were actually for instances with up to ten jobs.


\FloatBarrier
\vspace{-0.5cm}

%% file: 283_impact_of_the_brkga_components.tex
\subsection{\okC{Impact of the BRKGA components}}
\label{sec:impact_brkga_components}

\okC{In this section, we analyze the impact of the BRKGA-R-S-LS components.
We executed some additional experiments considering two versions of the BRKGA that remove different BRKGA-R-S-LS components.
The first, referred to as BRKGA with restarts and shakes (\mbox{BRKGA-R-S}), removes the local search component.
The second, denoted as BRKGA with restarts and local search (\mbox{BRKGA-R-LS}), removes the shaking perturbation component.
Moreover, the experiments were run on a small subset of the benchmark instances.
This subset contains nine randomly selected instances, one for each variation on the number of jobs ($n \in \{25, 50, 100\}$) and categories (S, M, L).
Furthermore, we performed ten independent executions of the BRKGA-R-S and the BRKGA-R-LS for each instance.
The results for the main approaches (BRKGA-R and BRKGA-R-S-LS) are the same as those reported in previous sections.}

\okC{As for the main BRKGA approaches (BRKGA-R and BRKGA-R-S-LS), we employed the \texttt{irace} package to define the parameters of the BRKGA-R-S and the BRKGA-R-LS.
Table \ref{tab:summary_parameters_brkga_components} displays the parameter settings suggested by \texttt{irace}.
Like the BRKGA-R-S-LS, we assume $r_{iLS} = n$ for the \mbox{BRKGA-R-LS}.
Besides, we also assume $\gamma^{reset} \in \{$OB, BI, BMS$\}$.
Finally, we highlight that setting $\gamma^{reset} =$ CB is equivalent to assuming $\gamma^{reset} =$ OB for the \mbox{BRKGA-R-LS} as well as for the BRKGA-R.
For additional details concerning the parameters and configurations, see Appendix~\ref{app:parameter_settings}.}

\begin{table}[!htp]\centering
\caption{\okC{Settings suggested by \texttt{irace} for the set of experiments on the BRKGA-R-S and the BRKGA-R-LS.}}
\label{tab:summary_parameters_brkga_components}
\scriptsize
\setstretch{0.85}
\okC{
\begin{tabular}{p{2cm}l c cc l cc}
\toprule
\multirow{2}{*}{Component} & \multirow{2}{*}{Notation} &  & \multicolumn{2}{c}{BRKGA-R-S} &  & \multicolumn{2}{c}{BRKGA-R-LS} \\ 
\cmidrule{4-5} \cmidrule{7-8} 
 &  &  & Used & \multicolumn{1}{c}{$\texttt{irace}$} &  & Used & \multicolumn{1}{c}{$\texttt{irace}$} \\
 \cmidrule{1-2} \cmidrule{4-5} \cmidrule{7-8} 
\multirow{5}{*}{Standard} & $n$ &  & $\checkmark$ & \multicolumn{1}{c}{--} &  & $\checkmark$ & -- \\
 & $p$ &  & $\checkmark$ & 187 &  & $\checkmark$ & 105 \\
 & $p_e$ &  & $\checkmark$ & 0.22 &  & $\checkmark$ & 0.35 \\
 & $p_m$ &  & $\checkmark$ & 0.10 &  & $\checkmark$ & 0.21 \\
 & $\rho_e$ &  & $\checkmark$ & 0.60 &  & $\checkmark$ & 0.68 \\
 \cmidrule{1-2} \cmidrule{4-5} \cmidrule{7-8} 
\multirow{3}{*}{$\parbox{3cm}{Initial solution\\generator}$} & $\alpha$ &  & $\checkmark$ & 0.01 &  & $\checkmark$ & 0.03 \\
 & $\lambda^{ws}$ &  & $\checkmark$ & 0.90 &  & $\checkmark$ & 0.58 \\
 & $n_{msi}$ &  & $\checkmark$ & 679 &  & $\checkmark$ & 313 \\
 \cmidrule{1-2} \cmidrule{4-5} \cmidrule{7-8} 
\multirow{8}{*}{Perturbation} & $n_{nimp}$ &  & -- & -- &  & $\checkmark$ & 144 \\
 & $R$ &  & $\checkmark$ & 141 &  & -- & -- \\
 & $R^{*}$ &  & $\checkmark$ & 5 &  & -- & -- \\
 & $R^{**}$ &  & $\checkmark$ & 9 &  & -- & -- \\
 & $s^{type}$ &  & $\checkmark$ & CHANGE &  & -- & -- \\
 & $\gamma^{weak}$ &  & $\checkmark$ & CB &  & -- & -- \\
 & $\gamma^{strong}$ &  & $\checkmark$ & OB &  & -- & -- \\
 & $\gamma^{reset}$ &  & $\checkmark$ & BMS &  & $\checkmark$ & OB \\ 
 \cmidrule{1-2} \cmidrule{4-5} \cmidrule{7-8} 
\multirow{4}{*}{Local search} & $b$ &  & -- & -- &  & $\checkmark$ & 7 \\
 & $L$ &  & -- & -- &  & $\checkmark$ & $0.18 \cdot n$ \\
 & $r_{pLS}$ &  & -- & -- &  & $\checkmark$ & 7 \\
 & $r_{iLS}$ &  & -- & -- &  & $\checkmark$ & -- \\ 
\bottomrule
\end{tabular}
}
\end{table}
\FloatBarrier

\begin{figure}[h]
    \centering
    \includegraphics[width=0.5\linewidth]{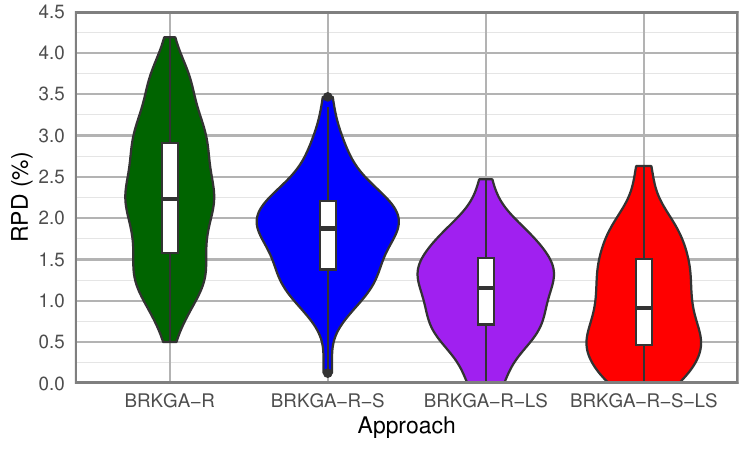}
    \caption{\okC{Impact of the BRKGA components.}}
    \label{fig:impact_brkga_components}
\end{figure}

\okC{Figure \ref{fig:impact_brkga_components} depicts the distribution of the RPD values for the nine instances selected, separated by each one of the four BRKGA approaches.
Each violin plot has a bandwidth equal to 0.25.
The RPD is calculated as: $100 \cdot \frac{obj^{run} - obj^{best}}{obj^{best}}$, where $obj^{run}$ is the makespan achieved on a run of the respective BRKGA and $obj^{best}$ is the best makespan acquired across the four approaches.
Each violin plot contains 90 observations, one for each of the ten independent runs of the respective BRKGA on each of the nine instances.
Notice that the medians decrease in the following order: BRKGA-R, BRKGA-R-S, BRKGA-R-LS, BRKGA-R-S-LS.
Besides, observe that both the shaking and the local search components improve the solutions of the BRKGA-R, especially the local search component.
Concerning the range of the RPD values, both the BRKGA-R-LS and the BRKGA-R-S-LS are generally competitive, although the bump arising in the lower tail of the BRKGA-R-S-LS distribution demonstrates the superiority of our main approach.
Thereupon, the BRKGA components, combined or not, positively contribute to the quality of the BRKGA solutions.}

\FloatBarrier

%% file: 300_finalcomments.tex
\section{Concluding remarks}
\label{sec:concluding}

In this paper, we considered the single-machine coupled task scheduling problem (SMCTSP) with exact delays to minimize the makespan.
We proposed a constraint programming (CP) model and a biased random-key genetic algorithm (BRKGA) to tackle the problem.
\okC{Our CP model applies well-established global constraints.
Conversely, our BRKGA combines some successful components in the literature: an initial solution generator, periodical restarts and shakes, and a local search algorithm.}
To the best of our knowledge, this is the first time such approaches have been proposed in the context of coupled task scheduling.

The computational experiments indicate that the proposed approaches are promising for tackling the problem. 
We summarize the results only for the 180 instances with at least 15 jobs, which are the ones that were not all solved to optimality by the approaches in the literature. 
\okC{Concerning the exact formulations, the experiments show that the \okC{CP model} outperformed the state-of-the-art MIP formulations, obtaining the smallest objective values, the smallest gaps, and the largest number of optimal solutions.
Moreover, only the \okC{CP}, compared to the MIP formulations, obtained feasible solutions for all the instances.}
Comparison between the tested approaches under the same computational settings, i.e., low time limit and one thread, indicated a superior performance of the BRKGA approach in terms of solution quality\okC{, compared to the CP model}.
\okC{Besides, the BRKGA approach acquired high-quality approximate solutions compared to the best known solutions within low computational times.
The CP model with multiple threads and a longer running time of 3600 seconds obtained the best-known solution in 163 out of 180 (90.56\%) of these instances.}
\okC{Finally, we show the importance of the shaking and local search components on the quality of our BRKGA solutions, especially the local search component.}

Our study also tackled an undesirable behavior that may be faced with a BRKGA for certain problems, as was the case with the SMCTSP, which is a fast population homogeneity.
We observed that some parameters of the BRKGA highly affect the population diversity through the generations.
Hence, we adopted a \textit{shaking} method.
Instead of resetting the whole population, which fully destroys useful parts of the random-key vector that usually appear in better solutions, it controls the intensity of the modifications in the elite set and only resets the non-elite set.
Combined with the injection of certain types of solution after each shake/reset, and a local search in selected situations, the resulting BRKGA significantly improved the overall solutions of the standard BRKGA with restarts.

%% file: 800_appendix_irace.tex
\section{Parameter settings}
\label{app:parameter_settings}

This section details our process of parameter selection for the BRKGA.
Recall that both BRKGA approaches have many parameters, especially BRKGA-R-S-LS.
Thus, instead of performing a full factorial design, we used the \texttt{irace} package \cite{iracepackage} version 4.2.0 implemented in R, which provides an automatic configuration tool for tuning optimization algorithms.
According to \citeA{iracepackage}, this software package implements the \textit{iterated racing}, a method that consists of three steps:
(1) sampling new configurations according to a particular distribution;
(2) selecting the best configurations from the newly sampled ones by means of racing; and 
(3) updating the sampling distribution to bias the sampling towards the best configurations. These three steps are repeated until a termination criterion is met.

For this task, we generated a new training set of instances for parametrization according to the general benchmark set proposed in \citeA{csbenchmark,KhaSalChe20}, using Python 3.11.9.
This new set contains 24 instances, two for each of the three categories and four different numbers of jobs, i.e., $n \in \{15, 25, 50, 100\}$.
In total, there are 24 instances.
We defined a block of instances of size four, which contains instances of the same category for each value of $n$.
This means that \texttt{irace} will only eliminate a configuration after evaluating a complete block and never in the middle of the block.
We also run the \texttt{irace} for each BRKGA approach, being that we used a budget of 4000 experiments for the BRKGA-R, \okC{5000 experiments for both BRKGA-R-S and BRKGA-R-LS,} and 6000 experiments for the BRKGA-R-S-LS, given that it contains more parameters.
The remaining \texttt{irace} parameters were configured according to the default settings.

Table \ref{tab:ranges_irace} exhibits the ranges used to tune the BRKGA parameters using \texttt{irace}.
Recall Sections \ref{sec:perturbation} and \ref{sec:tested_approaches_parameter_settings} for the meaning of the notations OB (overall best), CB (current best), \okC{BI (best initial solution)}, and BMS (best multi-start). 
Besides, we forbid $R^{**} < R^{*}$ and the maximum number of digits of the real numbers is two.
Finally, the parameter $L$ in the local search component was set to be proportional to $n$, i.e., $L = \lambda^{pLS} \cdot n$, where $\lambda^{pLS} \in [0,1]$ is the percentage of the number of jobs that defines the frequency of the periodic LS.
In this way, the greater the number of jobs, the lower the frequency of the periodic LS.
This aimed at reducing the total time spent in LS.

\begin{table}[h]
    \scriptsize
    \centering
    \caption{Ranges used to tune the BRKGA parameters using \texttt{irace}} \label{tab:ranges_irace}
    \resizebox{\textwidth}{!}{
    \okC{
    \begin{tabular}{lllclclclc}
        \toprule
        Component & \multicolumn{1}{c}{Notation} &  & BRKGA-R &  & BRKGA-R-S &  & BRKGA-R-LS &  & BRKGA-R-S-LS \\ \cmidrule{1-2} \cmidrule{4-4} \cmidrule{6-6} \cmidrule{8-8} \cmidrule{10-10} 
        \multirow{4}{*}{Standard} & $p$ &  & $\mathbb{Z} \cap [100,200]$ &  & $\mathbb{Z} \cap [100,200]$ &  & $\mathbb{Z} \cap [100,200]$ &  & $\mathbb{Z} \cap [100,200]$ \\
         & $p_e$ &  & $\mathbb{R} \cap [0.10, 0.50]$ &  & $\mathbb{R} \cap [0.10, 0.50]$ &  & $\mathbb{R} \cap [0.10, 0.50]$ &  & $\mathbb{R} \cap [0.10, 0.50]$ \\
         & $p_m$ &  & $\mathbb{R} \cap [0.10, 0.50]$ &  & $\mathbb{R} \cap [0.10, 0.50]$ &  & $\mathbb{R} \cap [0.10, 0.50]$ &  & $\mathbb{R} \cap [0.10, 0.50]$ \\
         & $\rho_e$ &  & $\mathbb{R} \cap [0.50, 0.80]$ &  & $\mathbb{R} \cap [0.50, 0.80]$ &  & $\mathbb{R} \cap [0.50, 0.80]$ &  & $\mathbb{R} \cap [0.50, 0.80]$ \\ 
         \cmidrule{1-2} \cmidrule{4-4} \cmidrule{6-6} \cmidrule{8-8} \cmidrule{10-10} 
        \multirow{3}{*}{\normalsize$\substack{\text{Initial solution} \\\text{generator}}$} & $\alpha$ &  & $\mathbb{R} \cap [0.00, 0.10]$ &  & $\mathbb{R} \cap [0.00, 0.10]$ &  & $\mathbb{R} \cap [0.00, 0.10]$ &  & $\mathbb{R} \cap [0.00, 0.10]$ \\
         & $\lambda^{ws}$ &  & $\mathbb{R} \cap [0.05, 1.00]$ &  & $\mathbb{R} \cap [0.05, 1.00]$ &  & $\mathbb{R} \cap [0.05, 1.00]$ &  & $\mathbb{R} \cap [0.05, 1.00]$ \\
         & $n_{msi}$ &  & $\mathbb{Z} \cap [300,1500]$ &  & $\mathbb{Z} \cap [300,1500]$ &  & $\mathbb{Z} \cap [300,1500]$ &  & $\mathbb{Z} \cap [300,1500]$ \\ \cmidrule{1-2} \cmidrule{4-4} \cmidrule{6-6} \cmidrule{8-8} \cmidrule{10-10} 
        \multirow{8}{*}{Perturbation} & $n_{nimp}$ &  & $\mathbb{Z} \cap [50,1000]$ &  & --- &  & $\mathbb{Z} \cap [50,1000]$ &  & - \\
         & $R$ &  & --- &  & $\mathbb{Z} \cap [50,250]$ &  & --- &  & $\mathbb{Z} \cap [50,250]$ \\
         & $R^{*}$ &  & --- &  & $\mathbb{Z} \cap [2,5]$ &  & --- &  & $\mathbb{Z} \cap [2,5]$ \\
         & $R^{**}$ &  & --- &  & $\mathbb{Z} \cap [3,10]$ &  & --- &  & $\mathbb{Z} \cap [3,10]$ \\
         & $s^{type}$ &  & --- &  & \{CHANGE, SWAP\} &  & --- &  & \{CHANGE, SWAP\} \\
         & $\gamma^{weak}$ &  & --- &  & \{CB, OB\} &  & --- &  & \{CB, OB\} \\
         & $\gamma^{strong}$ &  & --- &  & \{OB, BMS\} &  & --- &  & \{OB, BMS\} \\
         & $\gamma^{reset}$ &  & \{OB, BI, BMS\} &  & \{OB, BMS\} &  & \{OB, BI, BMS\} &  & \{OB, BMS\} \\ \cmidrule{1-2} \cmidrule{4-4} \cmidrule{6-6} \cmidrule{8-8} \cmidrule{10-10} 
        \multirow{3}{*}{Local search} & $b$ &  & --- &  & --- &  & $\mathbb{Z} \cap [1,10]$ &  & $\mathbb{Z} \cap [1,10]$ \\
         & $L$ &  & --- &  & --- &  & $\mathbb{R} \cap [0.10, 0.50] \cdot n$ &  & $\mathbb{R} \cap [0.10, 0.50] \cdot n$ \\
         & $r_{pLS}$ &  & --- &  & --- &  & $\mathbb{Z} \cap [2,7]$ &  & $\mathbb{Z} \cap [2,7]$\\
        \bottomrule
    \end{tabular}}}
\end{table}

\okC{Tables \ref{tab:irace_results_brkga_r}-\ref{tab:irace_results_brkga_r_s_ls}} exhibit the best configurations provided by \texttt{irace} for the BRKGA-R\okC{, BRKGA-R-S, BRKGA-R-LS,} and BRKGA-R-S-LS, listed from best to worst.

\begin{table}[!ht]
\scriptsize
\centering
\setlength{\tabcolsep}{4pt}
\caption{\texttt{irace} results for BRKGA-R, listed from best to worst}
\label{tab:irace_results_brkga_r}
\begin{tabular}{cccc l ccc l cc}
\toprule
\multicolumn{4}{c}{Standard} & & \multicolumn{3}{c}{\normalsize$\substack{\text{Initial solution} \\\text{generator}}$} & & \multicolumn{2}{c}{Perturbation} \\
\cmidrule{1-4} \cmidrule{6-8} \cmidrule{10-11}
$p$ & $p_e$ & $p_m$ & $\rho_e$ & & $\alpha$ & $\lambda^{ws}$ & $n_{msi}$ & & $n_{nimp}$ & $\gamma^{reset}$ \\
\cmidrule{1-4} \cmidrule{6-8} \cmidrule{10-11}
163 & 0.39 & 0.20 & 0.74 & & 0.02 & 0.28 & 1374 & & 956 & BMS \\
166 & 0.47 & 0.14 & 0.73 & & 0.03 & 0.53 & 1339 & & 859 & BMS \\
164 & 0.39 & 0.14 & 0.65 & & 0.01 & 0.57 & 1222 & & 568 & BMS \\
149 & 0.43 & 0.25 & 0.74 & & 0.01 & 0.35 & 1144 & & 879 & BMS \\
144 & 0.39 & 0.22 & 0.75 & & 0.03 & 0.23 & 1185 & & 974 & BMS \\
\bottomrule
\end{tabular}
\end{table}

\begin{table}[h]
\scriptsize
\centering
\setlength{\tabcolsep}{4pt}
\caption{\okC{\texttt{irace} results for BRKGA-R-S, listed from best to worst}}
\label{tab:irace_results_brkga_r_s}
\okC{
\begin{tabular}{cccc l ccc l ccccccc}
\toprule
\multicolumn{4}{c}{Standard} & & \multicolumn{3}{c}{\normalsize$\substack{\text{Initial solution} \\\text{generator}}$} & & \multicolumn{7}{c}{Perturbation} \\
\cmidrule{1-4} \cmidrule{6-8} \cmidrule{10-16}
$p$ & $p_e$ & $p_m$ & $\rho_e$ & & $\alpha$ & $\lambda^{ws}$ & $n_{msi}$ & & $R$ & $R^*$ & $R^{**}$ & $s^{type}$ & $\gamma^{weak}$ & $\gamma^{strong}$ & $\gamma^{reset}$ \\
\cmidrule{1-4} \cmidrule{6-8} \cmidrule{10-16}
187 & 0.22 & 0.10 & 0.60 & & 0.01 & 0.90 & 679 & & 141 & 5 & 9 & CHANGE & CB & OB & BMS \\
145 & 0.24 & 0.10 & 0.68 & & 0.02 & 0.90 & 891 & & 206 & 3 & 9 & CHANGE & CB & OB & BMS \\
\bottomrule
\end{tabular}}
\end{table}

\begin{table}[h]
\scriptsize
\centering
\setlength{\tabcolsep}{4pt}
\caption{\okC{\texttt{irace} results for BRKGA-R-LS, listed from best to worst}}
\label{tab:irace_results_brkga_r_ls}
\okC{
\begin{tabular}{cccc l ccc l cc l ccc}
\toprule
\multicolumn{4}{c}{Standard} & & \multicolumn{3}{c}{\normalsize$\substack{\text{Initial solution} \\\text{generator}}$} & & \multicolumn{2}{c}{Perturbation}  & & \multicolumn{3}{c}{Local search} \\ 
\cmidrule{1-4} \cmidrule{6-8} \cmidrule{10-11} \cmidrule{13-15}
$p$      & $p_e$ & $p_m$ & $\rho_e$ & & $\alpha$   & $\lambda^{ws}$ & $n_{msi}$ & & $n_{nimp}$   & $\gamma^{reset}$ & & $b$          & $L$  & $r_{pLS}$ \\ 
\cmidrule{1-4} \cmidrule{6-8} \cmidrule{10-11} \cmidrule{13-15}
105 & 0.35 & 0.21 & 0.68 & & 0.03 & 0.58 & 313 & & 144 & OB & & 7 & 0.18 & 7 \\
108 & 0.40 & 0.25 & 0.72 & & 0.04 & 0.62 & 433 & & 143 & OB & & 8 & 0.24 & 7 \\
182 & 0.30 & 0.19 & 0.69 & & 0.02 & 0.33 & 1044 & & 127 & OB & & 7 & 0.28 & 7 \\
195 & 0.43 & 0.17 & 0.72 & & 0.02 & 0.08 & 1475 & & 257 & OB & & 7 & 0.30 & 7 \\
\bottomrule
\end{tabular}}
\end{table}

\begin{table}[!ht]
\centering
\setlength{\tabcolsep}{4pt}
\scriptsize
\caption{\texttt{irace} results for BRKGA-R-S-LS, listed from best to worst}
\label{tab:irace_results_brkga_r_s_ls}
\begin{tabular}{cccc l ccc l ccccccc l ccc}
\toprule
\multicolumn{4}{c}{Standard} & & \multicolumn{3}{c}{\normalsize$\substack{\text{Initial solution} \\\text{generator}}$} & & \multicolumn{7}{c}{Perturbation} & & \multicolumn{3}{c}{Local search} \\
\cmidrule{1-4} \cmidrule{6-8} \cmidrule{10-16} \cmidrule{18-20}
$p$ & $p_e$ & $p_m$ & $\rho_e$ & & $\alpha$ & $\lambda^{ws}$ & $n_{msi}$ & & $R$ & $R^*$ & $R^{**}$ & $s^{type}$ & $\gamma^{weak}$ & $\gamma^{strong}$ & $\gamma^{reset}$ & & $b$ & $L$ & $r_{pLS}$ \\
\cmidrule{1-4} \cmidrule{6-8} \cmidrule{10-16} \cmidrule{18-20}
185 &  0.43 &  0.24 &  0.78 & &  0.01 &  0.22 & 602 & & 154 & 2 & 9 & SWAP & OB & OB & OB & & 9 & 0.21 & 7 \\
180 &  0.39 &  0.14 &  0.79 & &  0.04 &  0.88 & 1281 & & 121 & 4 & 10 & CHANGE & CB & BMS & OB & & 10 & 0.13 & 6 \\
164 &  0.43 &  0.11 &  0.70 & &  0.05 &  0.27 & 1156 & & 155 & 3 & 4 & SWAP & OB & OB & OB & & 9 & 0.24 & 7 \\
191 &  0.37 &  0.18 &  0.74 & &  0.02 &  0.13 & 649 & & 218 & 2 & 4 & SWAP & OB & OB & OB & & 9 & 0.19 & 7 \\
154 &  0.39 &  0.19 &  0.74 & &  0.02 &  0.57 & 1299 & & 98 & 5 & 8 & CHANGE & CB & BMS & OB & & 7 & 0.17 & 6 \\
197 &  0.45 &  0.14 &  0.77 & &  0.01 &  0.11 & 523 & & 159 & 2 & 6 & SWAP & OB & OB & OB & & 9 & 0.14 & 7 \\
\bottomrule
\end{tabular}
\end{table}
\FloatBarrier

%% file: 810_population_homogeneity.tex
\section{Population homogeneity in the BRKGA}
\label{app:pop_homogeneity_brkga}

\begin{figure}[!ht]
    \centering
    \subfigure[$p_e = 20\%$, $p_m = 10\%$ and $\rho_e = 70\%$\label{fig:effect_params_diversification_solution_improvement:01}]{
        \includegraphics[width=0.48\linewidth]{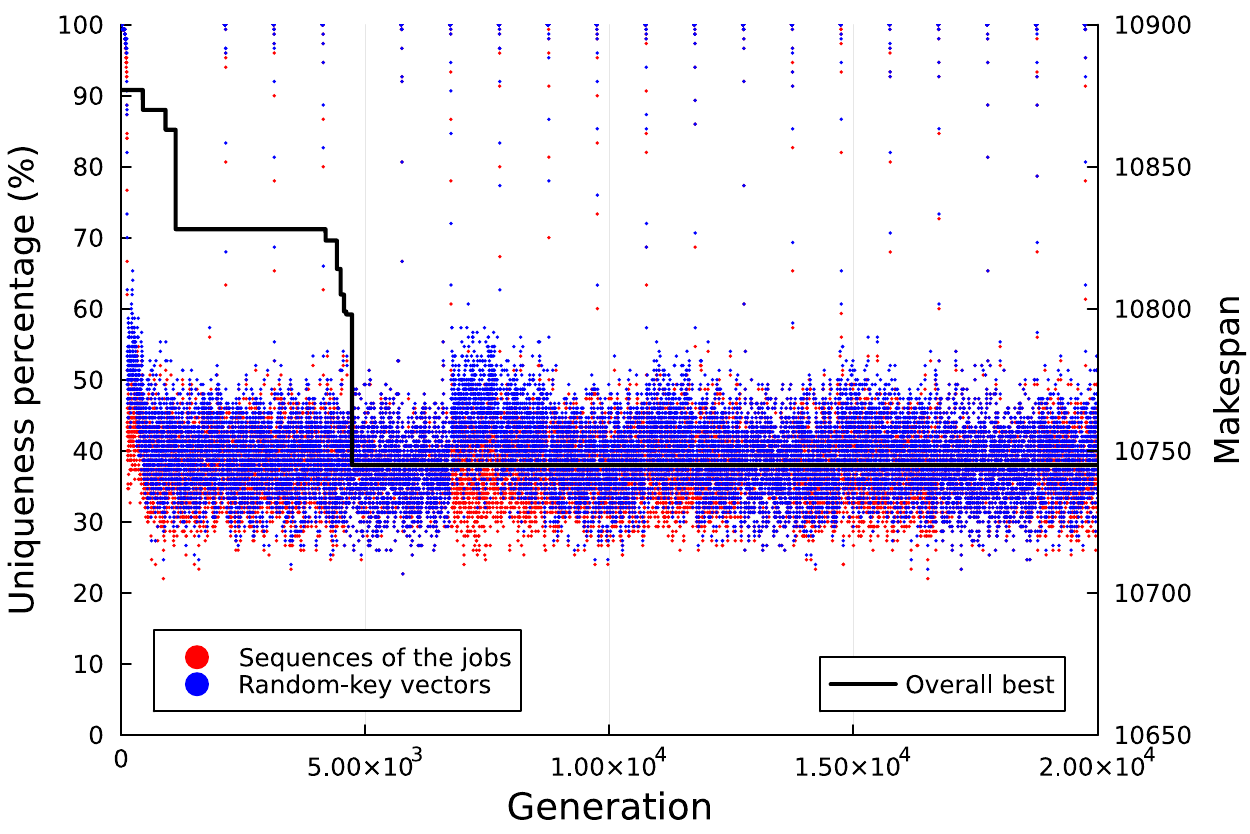}
    }
    \subfigure[$p_e = 10\%$, $p_m = 30\%$ and $\rho_e = 60\%$\label{fig:effect_params_diversification_solution_improvement:02}]{
        \includegraphics[width=0.48\linewidth]{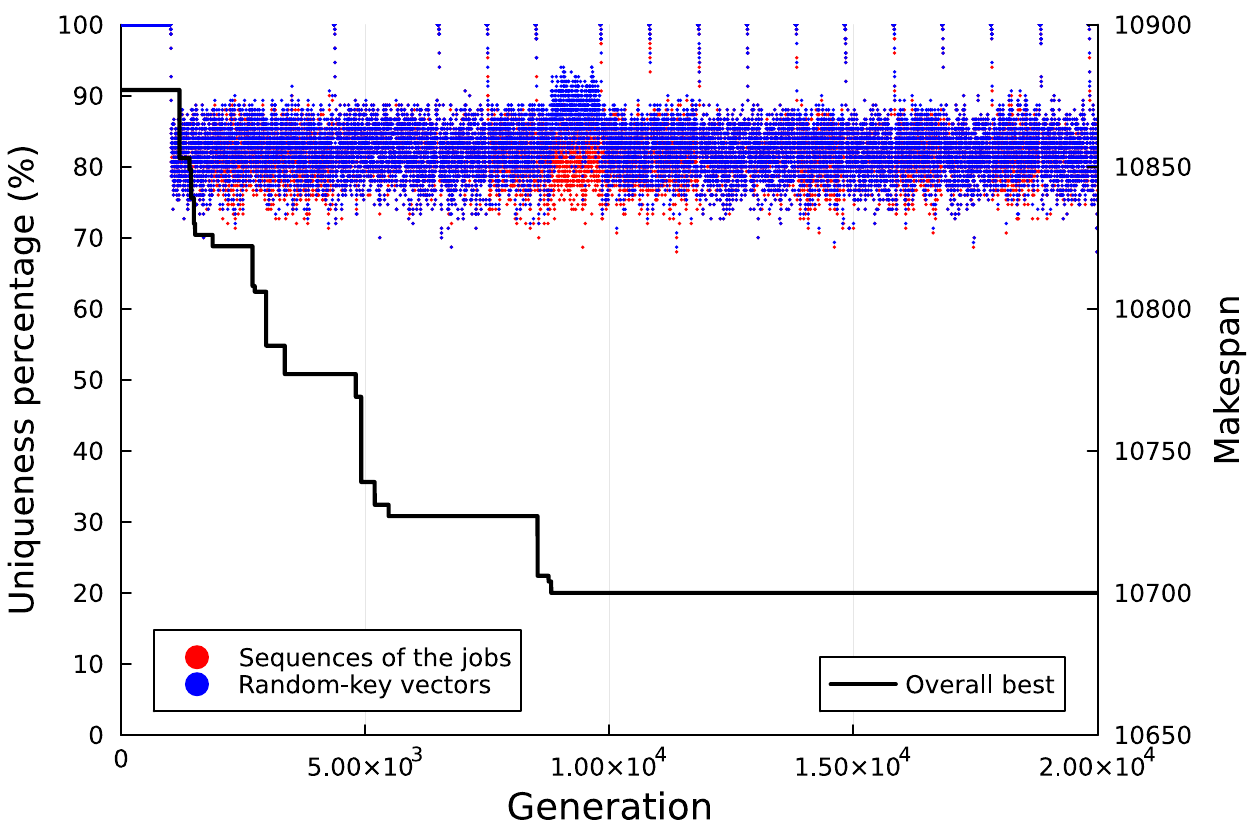}
    }
    \caption{Effect of $p_e$, $p_m$, and $\rho_e$ on BRKGA-R's population diversity and solution improvement. Instance 100\_2\_L\_gen.}
    \label{fig:effect_params_diversification_solution_improvement}
\end{figure}

In this section, we illustrate a fact observed during preliminary analysis regarding the percentage of unique random-key vectors in a population throughout the generations.
To accomplish that, we selected the instance 100\_2\_L\_gen from the general set of \citeA{csbenchmark, KhaSalChe20}'s benchmark.
Figure \ref{fig:effect_params_diversification_solution_improvement} depicts the effect of $p_e$, $p_m$, and $\rho_e$ on BRKGA-R's diversity and solution improvement.
For each population generation, the left axis shows the uniqueness percentage of the random-key vectors and the sequences of jobs.
Conversely, the right axis presents the evolution of the overall best solution in each generation.
We consider $p=150$, $n_{nimp} = 1000$, $\alpha = 0.05$, $\lambda^{ws} = 0.12$, $n_{msi} = 1000$, and $\gamma^{reset} =$ OB.

Observe that Figures \ref{fig:effect_params_diversification_solution_improvement:01} and \ref{fig:effect_params_diversification_solution_improvement:02} indicate that higher elite percentages $p_e$, lower mutant percentages $p_m$, and higher inheritance probability of elite parent key $\rho_e$ lead to a less diverse population.
Note that both configurations quickly converged to the same uniqueness percentage range after each reset.
Besides, in Figure \ref{fig:effect_params_diversification_solution_improvement:02} the uniqueness percentage did not converge before the first restart, which indicates that the percentage of warm-start solutions in the initial population also affected the convergence's velocity.
Furthermore, both plots indicate that several new series of improvement start a few generations after a new reset, particularly after a convergence.
Moreover, small key changes may not increase diversification as different random-key vectors convert to the same sequence when sorted.
Although the scenario presented in Figure \ref{fig:effect_params_diversification_solution_improvement:01} may be undesirable, it is not clear that the configuration of Figure \ref{fig:effect_params_diversification_solution_improvement:02} provides better results.

\begin{figure}[!ht]
    \centering
    \subfigure[$p_e = 20\%$, $p_m = 10\%$ and $\rho_e = 70\%$\label{fig:effect_params_diversification_elite_set:01}]{
        \includegraphics[width=0.48\linewidth]{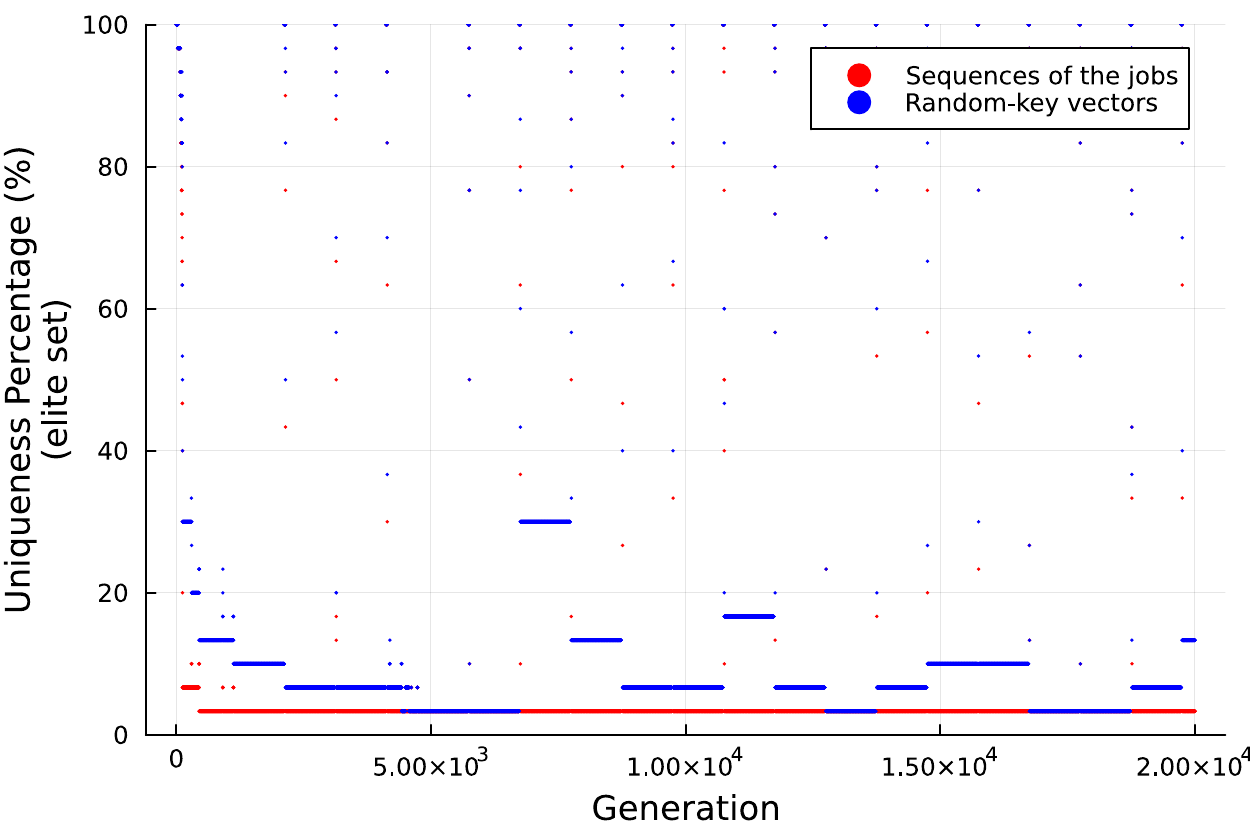}
    }
    \subfigure[$p_e = 10\%$, $p_m = 30\%$ and $\rho_e = 60\%$\label{fig:effect_params_diversification_elite_set:02}]{
        \includegraphics[width=0.48\linewidth]{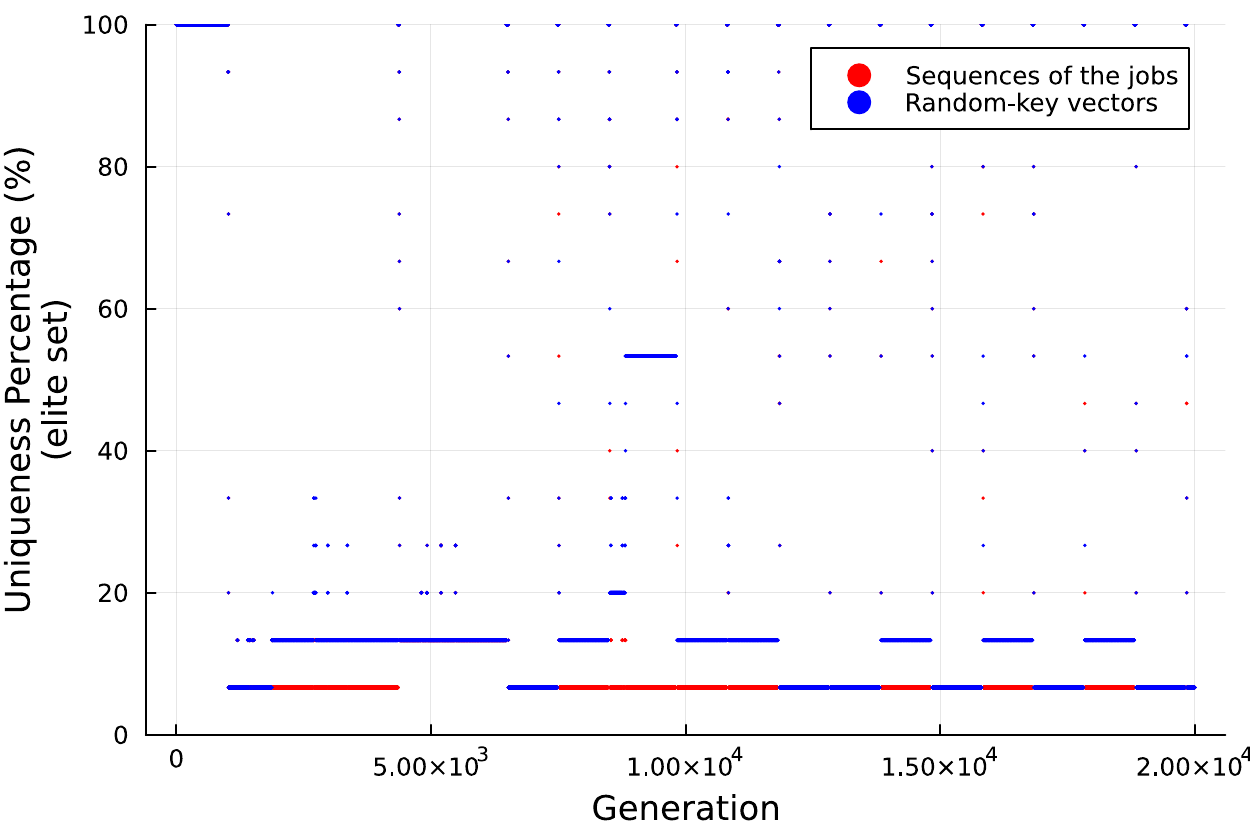}
    }
    \caption{Effect of $p_e$, $p_m$, and $\rho_e$ on BRKGA-R's elite set diversity. Instance 100\_2\_L\_gen.}
    \label{fig:effect_params_diversification_elite_set}
\end{figure}

In Figure \ref{fig:effect_params_diversification_solution_improvement:02}, we increased population diversity by setting a lower $p_e$, a higher $p_m$, and a lower $\rho_e$.
However, this figure does not depict whether these non-unique random-key vectors are concentrated in a few individuals or spread across the whole population.
Figure \ref{fig:effect_params_diversification_elite_set} deepens the analysis of the BRKGA-R's diversity in the elite set for the same runs of Figure \ref{fig:effect_params_diversification_solution_improvement}.
Observe that both configurations still resulted in homogeneous elite sets.
This situation is even more evident regarding the uniqueness percentage of the job's sequences.
Hence, modifying the parameters $p_e$, $p_m$, and $\rho_e$ may not be sufficient for generating a diverse population, particularly in the elite set.
